\pdfoutput=1 
\documentclass[OJMO,manuscript,Unicode]{cedram}
\newcommand{\R}{\mathbb{R}}
\newcommand{\tr}{\operatorname{tr}}
\newcommand{\rank}{\operatorname{rank}}
\newcommand{\supp}{\operatorname{supp}}
\newcommand{\conv}{\operatorname{conv}}
\makeatletter
\renewcommand*{\top}{%
  {\mathpalette\@transpose{}}%
}
\newcommand*{\@transpose}[2]{%
  \raisebox{\depth}{$\m@th#1\scriptscriptstyle\mathsf{T}$}%
}
\makeatother

\newcommand{\norm}[1]{\lVert #1\rVert}
\newcommand{\vol}{\varphi}
\usepackage{pgfplots}\usepgfplotslibrary{groupplots}\pgfplotsset{compat=1.17}

\title[Row-sparse generalized inverses with a few extra rows]{Row-sparse
  ah-symmetric reflexive generalized inverses with a few extra rows:
  local search with approximation guarantees}

\author{\firstname{Marcia} \lastname{Fampa}}
\address{Federal University of Rio de Janeiro\\ Brazil}
\email{fampa@cos.ufrj.br}

\author{\firstname{Jon} \lastname{Lee}}
\address{University of Michigan, Ann Arbor\\ USA}
\email{jonxlee@umich.edu}
\thanks{M.~Fampa was supported in part by CNPq grant 131219/2026-0.
  J.~Lee was supported in part by AFOSR grant FA9550-22-1-0172.}

\keywords{Moore--Penrose pseudoinverse, generalized inverse, row sparsity,
  2,1-norm, local search, approximation algorithm, least squares}

\begin{abstract}
  For a rank-$r$ matrix $A\in\R^{m\times n}$, ah-symmetric reflexive
  generalized inverses $H$ with few nonzero rows give least-squares solutions
  that are inexpensive to compute and use few variables. A known
  determinant-based local search finds one with exactly $r$ nonzero rows
  whose $2{,}1$-norm is within a factor $r$ of the minimum when the search
  stops at an exact local optimum, and this factor is essentially tight. We
  allow the nonzero rows of $H$ to lie in a set $K$ of $k\ge r$ indices,
  which also index columns of $A$, and study how to choose $K$. We bound the
  \emph{loss} of this restriction, the ratio of the restricted to the
  unrestricted minimum $2{,}1$-norm, by a simple geometric certificate, and
  we establish that a local search maximizing the squared rank-$r$ volume of
  $A[\,\cdot\,,K]$ gives loss at most $\sqrt{k(r+\epsilon k)/(k-r+1)}$,
  where $1+\epsilon$ is the swap tolerance. With $\epsilon=0$ this equals $r$
  at $k=r$ and is less than $\sqrt{2r}$ for $k\ge2r$; with $\epsilon=1/k$ the
  search runs in polynomial time on rational data. The same support also
  yields an ah-symmetric reflexive generalized inverse whose spectral norm
  is within a factor $\sqrt{1+(n-k)(r+\epsilon k)/(k-r+1)}$ of that of the
  Moore--Penrose pseudoinverse. The guarantees on the loss hold for every
  norm on the rows. Our bound on the certificate cannot be improved in
  general by more than a factor depending only on $k/r$, and on some
  instances, for a suitable row norm, neither can the bound on the loss; for
  the $2{,}1$-norm itself, however, the certificate can overestimate the
  loss by a factor $\sqrt r$. The local search is the Fedorov exchange
  algorithm of $0/1$ D-optimal design; we use this connection for a fast
  implementation and, through the Kiefer--Wolfowitz theorem, to establish that
  supports of at most $r(r+1)/2+1$ rows can have certificate at most $\sqrt r$, the
  best possible in general. On instances from earlier work, the observed
  loss is far below the guarantee, and the method matches two ADMM
  heuristics when few rows are removed and otherwise gives $2{,}1$-norms up
  to $27\%$ smaller, one to two orders of magnitude faster.
\end{abstract}

\def\currentvolume{0}

\makeatletter
\let\arxiv@copyrightline\copyrightline
\renewcommand\copyrightline{\setbox\z@\hbox{\arxiv@copyrightline}%
  \hbox to\wd\z@{\vrule\@width\z@\@height\ht\z@\@depth\dp\z@\hss}}
\let\arxiv@footnoterule\footnoterule
\renewcommand\footnoterule{\setbox\z@\vbox{\arxiv@footnoterule}\vbox to\ht\z@{}}
\makeatother

\begin{document}
\maketitle

\section{Introduction}\label{sec:intro}

For $A\in\R^{m\times n}$, Penrose \cite{Penrose1955} established that the
Moore--Penrose pseudoinverse $A^\dagger$ is the unique $H\in\R^{n\times m}$
satisfying
\[
  \text{(P1)}\ AHA=A,\qquad
  \text{(P2)}\ HAH=H,\qquad
  \text{(P3)}\ (AH)^\top=AH,\qquad
  \text{(P4)}\ (HA)^\top=HA.
\]
Following \cite{Rohde1964,XFLP2021}, $H$ is a \emph{generalized inverse} if it
satisfies (P1), it is \emph{reflexive} if it additionally satisfies (P2),
and it is \emph{ah-symmetric} if it satisfies (P3). Every generalized
inverse has $\rank H\ge\rank A$ (as $A=AHA$), and Rohde \cite{Rohde1964} established that
equality holds if and only if $H$ is reflexive: reflexive means of minimum
rank. If $H$ is an ah-symmetric generalized inverse, then
$\hat x:=Hb$ solves $\min\{\norm{Ax-b}_2: x\in\R^n\}$ for every $b$. When $H$
has few nonzero rows, the products $Hb$ are inexpensive to compute and the fitted linear
model $\hat x$ uses few of the $n$ variables (features). This motivates the line of work in
\cite{FampaLee2018,XFLP2021,FLPX2021,PFLX2024,PFLX2025} on \emph{row-sparse}
ah-symmetric reflexive generalized inverses with controlled norm.

Let $r:=\rank A$. Every generalized inverse has rank at least $r$, hence at
least $r$ nonzero rows, and one with exactly $r$ nonzero rows has rank
exactly $r$, hence is reflexive. So reflexivity is automatic for the
sparsest generalized inverses; for $k>r$ rows it is a genuine restriction,
built into the parametrization of Section~\ref{sec:red}. Two families of
methods are available.
\begin{itemize}
\item \emph{Minimum row count.} The column-block construction of
  \cite[Theorem~3.3]{XFLP2021} places $A[\,\cdot\,,T]^\dagger$ in the rows
  $T$, $|T|=r$. By \cite[Theorem~3.7]{XFLP2021} and
  \cite[Theorem~11]{PFLX2025}, choosing $T$ by a local search for
  $|\det A[S,T]|$, with tolerance $1+\epsilon$ per swap, gives a
  factor-$r(1+\epsilon)$ approximation for the $1$-norm and for the
  $2{,}1$-norm, respectively, and both factors are essentially best possible. Measured with tolerance $1+\epsilon$ on the
  squared determinant, as in Section~\ref{sec:vol}, the
  same factor reads $r\sqrt{1+\epsilon}$.
\item \emph{Norm minimization.} By \cite{PFLX2024}, minimizing
  $\norm{H}_{2{,}1}$ subject to (P1) is a convex problem whose optimal
  solutions automatically satisfy (P2) and (P3). The ADMM of \cite{PFLX2025}
  solves large instances efficiently, but its solutions typically have many more than $r$
  nonzero rows.
\end{itemize}
To bridge the two approaches, \cite{PFLX2025} proposed ADMM heuristics for the nonconvex
problem of minimizing (or merely achieving) a target number $k$ of nonzero
rows. These work well in practice but carry no guarantee. This paper builds
directly on \cite{PFLX2025}, and we assume familiarity with it: we use its
reduction of the $2{,}1$-norm problem, its ADMM algorithms, and its test
instances.

A third direction, taken in \cite{FLP2021}, keeps the norm small and relaxes
reflexivity instead. Starting from a $1$-norm minimizing ah-symmetric
generalized inverse, which need not be reflexive and typically has rank
larger than $r$, \cite{FLP2021} imposes reflexivity gradually and so trades
$1$-norm and sparsity against rank, ending at rank $r$. The present paper is
complementary: we keep reflexivity, that is, minimum rank $r$, and relax
instead the number of nonzero rows from $r$ to $k$.

\medskip\noindent\textbf{Contribution.} We study the problem of choosing a row support $K$
of size $k$, with $r\le k\le n$, and then minimizing the $2{,}1$-norm over
ah-symmetric reflexive generalized inverses supported on $K$. Our results are
as follows.
\begin{enumerate}
\item A \emph{certificate} (Theorem~\ref{thm:cert}): the \emph{loss} of
  restricting to $K$, that is, the ratio $z(K)/z^*$ of the restricted to the
  unrestricted minimum $2{,}1$-norm (see Section~\ref{sec:red}), is at most the largest, over the columns $a_l$ with $l\notin K$, of the
  smallest $\ell_1$-norm of a vector $c$ with $A[\,\cdot\,,K]c=a_l\,$, that
  is, of a coefficient vector that represents $a_l$ in terms of the
  columns of $A$ indexed by $K$. This
  quantity depends only on $K$, and the argument of \cite[Theorem~11]{PFLX2025}
  is its special case $k=r$.
\item A \emph{guarantee} (Theorem~\ref{thm:main}): a local maximizer of the
  squared rank-$r$ volume over $k$-subsets has loss at most
  $\sqrt{k(r+\epsilon k)/(k-r+1)}$. At $k=r$ it equals $r\sqrt{1+\epsilon}$. With
  $\epsilon=0$ it is at most $\sqrt{\kappa r/(\kappa-1)}$ for $k=\kappa r$ (here and below,
  $\kappa:=k/r$ is a ratio, not necessarily an integer), and
  at most $\sqrt{2r}$ for every $k\ge2r$. For rational $A$ and $\epsilon=1/k$, the
  local search runs in polynomial time (Proposition~\ref{prop:poly}).
\item A \emph{limitation} (Propositions~\ref{prop:coh}, \ref{prop:lower},
  and~\ref{prop:lowersharp}): if the columns of $A$ are unit vectors with
  mutual coherence $\operatorname{coh}(A)$ (the largest $|\langle a_i\,,a_j\rangle|$
  over distinct columns), every $k$-subset of rank $r$ with $k<n$ has
  certificate value at least $1/\operatorname{coh}(A)$, and unions of mutually unbiased bases give $1/\operatorname{coh}(A)=\sqrt r$ with
  only polynomially many columns. For general $r$, there are instances on
  which every $k$-subset
  has certificate value at least $\theta\sqrt{r/(2\ln(2k))}$, and also at
  least $\theta\sqrt{(1-\delta)r/(2\ln(5k/(r\delta^2)))}$, for any
  $\delta\in(0,1)$ and any $\theta<1$ (the instances depend on $\theta$),
  under explicit conditions on $k$ and $r$. The second bound, obtained
  from \v{S}id\'ak's inequality instead of a union bound (and sharpened
  further with the Mills ratio), has no logarithm in $r$.
  For $k=\kappa r$ the analysis behind Theorem~\ref{thm:main}
  is therefore tight up to an explicit factor depending only on $\kappa$:
  at most $4.5$ for $\kappa=2$ and every $r\ge7$, and at most
  $\sqrt{\kappa/(\kappa-1)}$ when $r$ is a power of $2$ and $\kappa<2$, or a
  power of $4$ and $\kappa<r/2+1$. For every $r\ge2$ and every $k$ with
  $r<k<2r$, an instance with $n=2r$ columns built from the discrete cosine
  transform gives the factor $\sqrt{2k/(k-r+1)}$, without requiring
  a Hadamard matrix (Example~\ref{ex:dct}). Conversely, a support of at most
  $r(r+1)/2+1$ rows taken from a D-optimal design has certificate at most
  $\sqrt r$ (Theorem~\ref{thm:dsupp}), so for supports of that size $\sqrt r$
  is exactly the best-possible bound.
  Theorem~\ref{thm:bases}
  separates the two notions: on unions of orthonormal bases, every support
  that contains one of the bases has $z(K)=z^*$, although, when the bases
  are mutually unbiased, its certificate is $\sqrt r$ unless $K$ is the
  full support;
  yet, on
  unions of $p\ge2$ mutually unbiased bases, global volume maximizers with $k=qr$
  rows, $1\le q<p$, lose exactly $\sqrt r=\rho(K)$ for a suitable norm on the
  rows. So the certificate is attained by the loss for some row norm, and,
  as a guarantee for all row norms (item~\ref{it:norms}), Theorem~\ref{thm:main} is
  tight for the loss itself up to the factor $\sqrt{k/(k-r+1)}$ when
  $r$ is a power of $4$ and $k$ is a multiple of $r$ with $2r\le k\le r^2/2$
  (Corollary~\ref{cor:bases}).
\item\label{it:norms} \emph{Other norms} (Theorem~\ref{thm:norms}): the certificate and the
  guarantee hold for $\sum_i\mu(H_{i\cdot})$ with any norm $\mu$, including
  the $1$-norm, for which the restricted problem is a linear program.
\item \emph{Experiments} (Section~\ref{sec:exp}), on the instances of
  \cite{PFLX2025}, demonstrating that the observed loss can be far below both the guarantee
  and the certificate. Compared with the two ADMM heuristics of
  \cite{PFLX2025} that target the same number of rows, the method, with a
  simple polishing step, is comparable to the better one when few
  rows are removed, and gives $2{,}1$-norms $2.5\%$ to $27\%$ smaller, $68$ to
  $394$ times faster on the same machine, when the target is close to $r$.
\end{enumerate}
We also establish that determinant local search does not depend on the
choice of row set $S$ (Proposition~\ref{prop:det}), and that, for $k<n$, a volume local
maximizer $K$ also yields an ah-symmetric reflexive generalized inverse
supported on $K$, namely the one of minimum Frobenius norm, whose spectral
norm is within a factor $\sqrt{1+(n-k)(r+\epsilon k)/(k-r+1)}$ of
$\norm{A^\dagger}_2$ (Proposition~\ref{prop:spectral}); at $k=r$ this
recovers a bound that \cite{Pan2000} implies.

\medskip\noindent\textbf{What is inherited and what is new.} Several
ingredients are known. \cite{PFLX2024} gave the parametrization of
ah-symmetric reflexive generalized inverses (Lemma~\ref{lem:param}). The
swap formula (Lemma~\ref{lem:swap}) is standard in D-optimal design, and
\cite{PFL2026} uses it for fast local search. For $\epsilon=0$, the bound on
leverage scores (Lemma~\ref{lem:lev}) is \cite[Lemma~9]{MSTX2019}, and we add
the tolerance $\epsilon$. \cite[Theorem~11]{PFLX2025} proved the
case $k=r$ of our guarantee (Corollary~\ref{cor:kr}). Berg \cite{Berg1986} is
the source of the convex-combination identity in Remark~\ref{rem:vs}, and
the scheme of the polynomial-time argument follows \cite{XFLP2021} closely.
Finally, \cite{PFL2026} and \cite{PFLX2025} supply, respectively, the
initialization and the two ADMM heuristics used in Section~\ref{sec:exp}.

Our new contributions are the following. Theorem~\ref{thm:cert} isolates
the dual-feasibility argument of \cite[Theorem~3.7]{XFLP2021} and
\cite[Theorem~11]{PFLX2025}, which was tied to $|K|=r$ and Cramer's rule, as
a certificate valid for every support $K$; this is what allows $k>r$.
Theorem~\ref{thm:main} extends the factor-$r$ guarantee to every $k\ge r$,
with a factor that drops to order $\sqrt r$. Proposition~\ref{prop:spectral}
extends a spectral-norm bound that \cite{Pan2000} implies from $k=r$ to every
$k\ge r$, improving it by roughly $\sqrt{k-r+1}$, and it adds an
a~posteriori bound. Proposition~\ref{prop:poly} bounds the number of swaps by
a quantity that depends on $k$ only logarithmically\ for fixed $\epsilon$. The lower bounds of
Section~\ref{sec:limits} are new: the earlier tightness results,
\cite[Example~A.1]{XFLP2021} and \cite[Remark~15]{PFLX2025}, concern only
$k=r$. Theorem~\ref{thm:dsupp} uses the Kiefer--Wolfowitz theorem to establish
that the coherence lower bound $\sqrt r$ is attained by supports of
at most $r(r+1)/2+1$ rows. Theorem~\ref{thm:norms} extends the $1$-norm guarantee of
\cite[Theorem~3.7]{XFLP2021} and the $2{,}1$-norm guarantee of
\cite[Theorem~11]{PFLX2025} to all row norms and all $k\ge r$; for the
$1$-norm at $k=r$, however, it measures against a weaker benchmark than the
one in \cite{XFLP2021} (Remark~\ref{rem:onenorm}).
Theorem~\ref{thm:bases} is also new: it establishes that the lower bounds
on the certificate do not transfer to the $2{,}1$-norm loss, and that they
do transfer, with equality, to the loss for a suitable row norm.

\medskip\noindent\textbf{Notation.} For a matrix $W$, $W_{i\cdot}$ is its $i$-th row,
$W[S,T]$ the submatrix with rows $S$ and columns $T$ (a dot means all),
$\norm{W}_{2{,}1}:=\sum_i\norm{W_{i\cdot}}_2\,$, and $\norm{W}_{2{,}0}$ the number
of nonzero rows of $W$. For $R\in\R^{r\times n}$, we write
$r_j:=R[\,\cdot\,,j]$ and $R_K:=R[\,\cdot\,,K]$. Finally, $M:=\{1,\ldots,m\}$ and
$N:=\{1,\ldots,n\}$, and a matrix with $n$ rows is \emph{supported on}
$K\subseteq N$ if its rows outside $K$ are zero. For $K\subseteq N$, $j\in K$,
and $l\notin K$, we write $K+l:=K\cup\{l\}$ and $K-j+l:=(K\setminus\{j\})\cup\{l\}$;
passing from $K$ to $K-j+l$ is a \emph{swap}.

\section{Reduction and the support-restricted problem}\label{sec:red}

Throughout, $A\ne0$, so $r\ge1$. Let $A=U_1DV_1^\top$ be a compact singular-value decomposition, with
$U_1\in\R^{m\times r}$ and $V_1\in\R^{n\times r}$ having orthonormal columns
and $D$ diagonal and positive, and define $R:=DV_1^\top\in\R^{r\times n}$, so
that $A=U_1R$.

\begin{lemm}[{cf.\ \cite[Lemmas~2--4]{PFLX2024}}]\label{lem:param}
  $H\in\R^{n\times m}$ is an ah-symmetric reflexive generalized inverse of
  $A$ if and only if $H=XU_1^\top$ for some $X\in\R^{n\times r}$ with $RX=I_r\,$.
  In that case $\rank H=r$, and $\norm{H_{i\cdot}}_2=\norm{X_{i\cdot}}_2$ for all
  $i$; in particular $H$ and $X$ have the same nonzero rows and the same
  $2{,}1$-norm.
\end{lemm}

\begin{proof}
  If $H=XU_1^\top$ with $RX=I$, then $AHA=U_1RXU_1^\top U_1R=U_1R=A$,
  $AH=U_1RXU_1^\top=U_1U_1^\top$ is symmetric, and $HAH=XU_1^\top U_1RXU_1^\top=H$.
  Conversely, by \cite[Lemmas~2--4]{PFLX2024}, every such $H$ has the form
  $(V_1D^{-1}+V_2Z)U_1^\top\,$, where $V_2$ spans the orthogonal complement of
  the range of $V_1$; then $X:=V_1D^{-1}+V_2Z$ satisfies $RX=I$. The norm
  identity holds because $U_1$ has orthonormal columns, and
  $\rank H=\rank X=r$ because $RX=I_r\,$.
\end{proof}

By \cite[Corollary~9]{PFLX2024}, every $2{,}1$-norm minimizing generalized
inverse is reflexive and ah-symmetric. With Lemma~\ref{lem:param}, this gives
\begin{equation}\label{eq:zstar}
  z^*:=\min\{\norm{H}_{2{,}1}: AHA=A\}
  =\min\{\norm{X}_{2{,}1}: RX=I_r\}.
\end{equation}
For $K\subseteq N$ with $|K|=k$ and $\rank R_K=r$, the
\emph{support-restricted problem} is
\begin{equation}\label{eq:zK}
  z(K):=\min\{\norm{X_K}_{2{,}1}: R_KX_K=I_r\,,\ X_K\in\R^{k\times r}\}.
\end{equation}
Padding an optimal $X_K$ with zero rows and multiplying by $U_1^\top$ gives an
ah-symmetric reflexive generalized inverse with at most $k$ nonzero rows and
$2{,}1$-norm $z(K)\ge z^*$. For $k=r$, problem~\eqref{eq:zK} has the unique
feasible point $X_K=R_K^{-1}\,$. Because $A[\,\cdot\,,K]=U_1R_K\,$, this is
exactly the column-block solution $A[\,\cdot\,,K]^\dagger$ of
\cite[Theorem~3.3]{XFLP2021}; for $|T|=r$ with $R_T$ nonsingular, we write
$H_T$ for the column-block solution over $T$, which has $A[\,\cdot\,,T]^\dagger$
in the rows $T$ and zeros elsewhere. For $k>r$, \eqref{eq:zK} is a small
second-order-cone problem with $k\times r$ unknowns, which the ADMM of
\cite[Algorithm~2]{PFLX2025} solves directly. To do so, replace $V_1D^{-1}$
by $R_K^\dagger$ and $V_2$ by an orthonormal basis of the null space of $R_K\,$.

\begin{prop}\label{prop:det}
  Let $S\subseteq M$ with $|S|=r$ index linearly independent rows of $A$.
  Then $U_1[S,\cdot\,]$ is nonsingular, and
  $\det(A[S,K]A[S,K]^\top)=\det(U_1[S,\cdot\,])^2\det(R_KR_K^\top)$ for every
  $K\subseteq N$.
  In particular, for $|T|=r$, $|\det A[S,T]|=|\det U_1[S,\cdot\,]|\cdot|\det R_T|$.
  Hence the local maximizers of $|\det A[S,T]|$ over swaps in $T$ used in
  \cite{XFLP2021,PFLX2025} do not depend on the choice of $S$; they are the
  local maximizers of $|\det R_T|$, equivalently of $|\det V_1[T,\cdot\,]|$.
\end{prop}

\begin{proof}
  $A[S,\cdot\,]=U_1[S,\cdot\,]R$ has rank $r$, so $U_1[S,\cdot\,]$ is
  nonsingular. Moreover $A[S,K]=U_1[S,\cdot\,]R_K\,$, and the identity
  follows from multiplicativity of the determinant.
\end{proof}

Proposition~\ref{prop:det} also establishes that determinant-based local search is
the ``local maximum volume'' pivoting familiar from the rank-revealing
factorizations of \cite{GTZ1997,Pan2000}; Proposition~\ref{prop:spectral}
below gives a corresponding spectral-norm guarantee.

\section{A certificate for the loss of restricting the support}\label{sec:cert}

By standard conic duality (the constraints are affine and~\eqref{eq:zK} is
feasible), we have
\begin{align}
  z^* &= \max\big\{\tr\Lambda:\ \Lambda\in\R^{r\times r},\
        \norm{\Lambda^\top r_j}_2\le 1,\ j\in N\big\},\label{eq:dual}\\
  z(K)&= \max\big\{\tr\Lambda:\ \Lambda\in\R^{r\times r},\
        \norm{\Lambda^\top r_j}_2\le 1,\ j\in K\big\},\label{eq:dualK}
\end{align}
with the maxima attained. The restricted dual~\eqref{eq:dualK} differs from
\eqref{eq:dual} only by dropping the constraints for $j\notin K$.

\begin{theo}\label{thm:cert}
  For $K$ with $\rank R_K=r$ and $l\notin K$, let
  $\beta_K(l):=\min\{\norm{c}_1: R_Kc=r_l\}$, and let
  $\rho(K):=\max\{1,\max_{l\notin K}\beta_K(l)\}$. Then
  $z^*\le z(K)\le\rho(K)\,z^*$.
  Equivalently, $\rho(K)$ is the least $\rho\ge1$ such that every column of
  $R$ lies in $\rho\cdot\conv\{\pm r_i: i\in K\}$.
\end{theo}

\begin{proof}
  Let $\Lambda$ be optimal for~\eqref{eq:dualK}. For $l\notin K$, choose
  $c$ with $R_Kc=r_l$ and $\norm{c}_1=\beta_K(l)$. Then
  \[
    \norm{\Lambda^\top r_l}_2=\Big\|\sum_{i\in K}c_i\Lambda^\top r_i\Big\|_2
    \le\sum_{i\in K}|c_i|\,\norm{\Lambda^\top r_i}_2\le\beta_K(l).
  \]
  Thus $\Lambda/\rho(K)$ is feasible for~\eqref{eq:dual}, and
  $z^*\ge\tr\Lambda/\rho(K)=z(K)/\rho(K)$.
\end{proof}

\begin{coro}[{\cite[Theorem~11]{PFLX2025}}]\label{cor:kr}
  If $|T|=r$, $R_T$ is nonsingular, and $|\det R_T|$ cannot be increased by more than a factor
  $1+\epsilon$ by a single swap, then $\norm{H_T}_{2{,}1}=z(T)\le r(1+\epsilon)z^*$.
\end{coro}

\begin{proof}
  As in Knuth's analysis \cite{Knuth1985}, with Cramer's rule we can establish that, for
  $l\notin T$, the unique $c$ with $R_Tc=r_l$ has
  $|c_i|=|\det R_{T-t_i+l}|/|\det R_T|\le1+\epsilon$, where $t_i$ is the
  $i$-th element of $T$. Hence
  $\beta_T(l)\le r(1+\epsilon)$, and Theorem~\ref{thm:cert} applies.
\end{proof}

Theorem~\ref{thm:cert} also suggests an inexpensive a~posteriori certificate. It
uses any approximately optimal dual for the small problem~\eqref{eq:zK}, as
produced, for example, by the scaled multiplier of ADMM.

\begin{prop}\label{prop:post}
  Let $\Theta\in\R^{k\times r}$ have rows of Euclidean norm at most one and
  columns in the range of $R_K^\top\,$. Set $\Lambda:=(R_KR_K^\top)^{-1}R_K\Theta$
  and $\nu:=\max\{1,\max_{j\in N}\norm{\Lambda^\top r_j}_2\}$. Then
  $z^*\ge\tr(\Lambda)/\nu$. If $\Theta=R_K^\top\Lambda^\star$ for an optimal
  solution $\Lambda^\star$ of~\eqref{eq:dualK}, then $\Lambda=\Lambda^\star$,
  $\tr\Lambda=z(K)$, and the certified ratio is $z(K)/z^*\le\nu\le\rho(K)$.
\end{prop}

\begin{proof}
  $R_K^\top\Lambda=\Theta$, so $\Lambda$ is feasible for~\eqref{eq:dualK}.
  Then $\Lambda/\nu$ is feasible for~\eqref{eq:dual}. The inequality
  $\nu\le\rho(K)$ is the proof of Theorem~\ref{thm:cert}.
\end{proof}

\section{Volume local search over \texorpdfstring{$k$}{k}-subsets}\label{sec:vol}

For $K\subseteq N$, define the squared rank-$r$ volume
\[
  \vol(K):=\det(R_KR_K^\top)=\sum_{T\subseteq K,\,|T|=r}\det(R_T)^2,
\]
where the equality is the Cauchy--Binet formula. Thus, when $\rank R_K=r$, $\vol(K)$ is the
product of the squared nonzero singular values of $A[\,\cdot\,,K]=U_1R_K\,$;
otherwise $\vol(K)=0$.
By Proposition~\ref{prop:det}, $\vol(K)$ is proportional to
$\det(A[S,K]A[S,K]^\top)$ for any fixed set $S$ of $r$ independent rows, so it
can be evaluated on rational data without an SVD.

In the language of optimal experimental design, maximizing $\vol$ over
$k$-subsets is $0/1$ D-optimal design: choose $k$ of $n$ candidate
experiments (\emph{design points}, here the columns of $R$) so as to
maximize the determinant of the resulting \emph{information matrix}, here
$\sum_{j\in K}r_jr_j^\top=R_KR_K^\top\,$. The local
search below is the classical swap-based heuristic for this problem, the
Fedorov exchange algorithm; see
Remark~\ref{rem:dopt} for the precise correspondence.

\begin{defi}\label{def:loc}
  Let $r\le k\le n$ and $\epsilon\ge0$. A $k$-subset $K$ with $\vol(K)>0$ is
  a \emph{$(1+\epsilon)$-local maximizer} if
  $\vol(K-j+l)\le(1+\epsilon)\vol(K)$ for all $j\in K$ and $l\notin K$.
\end{defi}

For a $k$-subset $K$ with $\vol(K)>0$, let $G:=R_KR_K^\top$ (nonsingular),
$\alpha_l:=r_l^\top G^{-1}r_l$ for $l\notin K$, $h_j:=r_j^\top G^{-1}r_j$ for
$j\in K$, and $p_{jl}:=r_j^\top G^{-1}r_l\,$.
Here $h_j$ is the \emph{leverage score} of $j$ in $K$, in the sense of the
regression diagnostics of \cite{HoaglinWelsch1978} and of the randomized
numerical linear algebra of \cite{DMMW2012}: the $j$-th diagonal
entry of the orthogonal projector $R_K^\top G^{-1}R_K$ onto the row space of
$R_K\,$, so $0\le h_j\le1$ and $\sum_{j\in K}h_j=r$. For $l\notin K$, $\alpha_l$
is the analogous quantity for a column outside $K$, and it can exceed $1$.
For $K=N$ we have $G=RR^\top=D^2$, and the leverage score of column $j$ of
$A$ is $r_j^\top(RR^\top)^{-1}r_j=\norm{V_1[j,\cdot\,]}_2^2\,$.

\begin{lemm}[{standard; see \cite[Section~4.4]{PFL2026}}]\label{lem:swap}
  For $j\in K$ and $l\notin K$,
  $\vol(K-j+l)=\vol(K)\,[(1+\alpha_l)(1-h_j)+p_{jl}^2]$.
\end{lemm}

\begin{proof}
  $R_{K-j+l}R_{K-j+l}^\top=G+r_lr_l^\top-r_jr_j^\top\,$. By the matrix determinant
  lemma,
  \[
    \det(G+r_lr_l^\top-r_jr_j^\top)
    =\det G\cdot\det\begin{pmatrix}1+\alpha_l & p_{jl}\\ -p_{jl} & 1-h_j\end{pmatrix}.
    \qedhere
  \]
\end{proof}

Lemma~\ref{lem:swap} lets all $k(n-k)$ swap ratios be evaluated from
$G^{-1}$ with matrix products, and $G^{-1}$ is maintained by rank-one
updates.

\begin{lemm}[{cf.\ \cite[Lemma~9]{MSTX2019}}]\label{lem:lev}
  If $K$ is a $(1+\epsilon)$-local maximizer, then
  $\alpha_l\le(r+\epsilon k)/(k-r+1)$ for every $l\notin K$.
\end{lemm}

\begin{proof}
  Fix $l\notin K$ and let $G':=G+r_lr_l^\top\,$. For $i\in K+l$, the leverage
  $h'_i:=r_i^\top G'^{-1}r_i$ satisfies $\sum_{i\in K+l}h'_i=\tr(G'^{-1}G')=r$.
  By the Sherman--Morrison formula, $h'_l=\alpha_l/(1+\alpha_l)$. Hence some $j\in K$
  has $h'_j\le(r-\alpha_l/(1+\alpha_l))/k$. Now $\vol(K+l)=(1+\alpha_l)\vol(K)$ and
  $\vol(K+l-j)=(1-h'_j)\vol(K+l)$, so local maximality gives
  $(1+\alpha_l)(1-h'_j)\le1+\epsilon$. Therefore
  \[
    (1+\alpha_l)-\frac{r(1+\alpha_l)-\alpha_l}{k}\le1+\epsilon,
    \quad\text{that is,}\quad
    \alpha_l\,\frac{k-r+1}{k}\le\epsilon+\frac rk .
    \qedhere
  \]
\end{proof}

\begin{theo}\label{thm:main}
  Let $r\le k<n$ and $\epsilon\ge0$, and let $K$ be a $(1+\epsilon)$-local
  maximizer of $\vol$ over $k$-subsets. Then every optimal solution
  of~\eqref{eq:zK}, padded and multiplied by $U_1^\top\,$, is an ah-symmetric
  reflexive generalized inverse $H$ of $A$ with at most $k$ nonzero rows and
  $\norm{H}_{2{,}1}=z(K)\le\sqrt{k(r+\epsilon k)/(k-r+1)}\;z^*$.
\end{theo}

\begin{proof}
  For $l\notin K$, let $c:=R_K^\dagger r_l=R_K^\top G^{-1}r_l\,$. Then $R_Kc=r_l$
  and $\norm{c}_2^2=r_l^\top G^{-1}r_l=\alpha_l\,$. Hence
  $\beta_K(l)\le\norm{c}_1\le\sqrt k\,\norm{c}_2\le\sqrt{k(r+\epsilon k)/(k-r+1)}$
  by Lemma~\ref{lem:lev}. This bound is at least $1$, so
  Theorem~\ref{thm:cert} applies.
\end{proof}

\begin{rema}\label{rem:values}
  At $k=r$, the factor is $r\sqrt{1+\epsilon}$. This agrees with
  Corollary~\ref{cor:kr} and \cite[Theorem~11]{PFLX2025}, which state the
  factor $r(1+\epsilon')$ for a local search on $|\det R_T|$ with tolerance
  $1+\epsilon'$: because $\vol(T)=\det(R_T)^2$, tolerance $1+\epsilon$ on
  $\vol$ is tolerance $\sqrt{1+\epsilon}$ on $|\det R_T|$. With $\epsilon=0$ and $k=\kappa r$,
  the factor is $\sqrt{\kappa r^2/((\kappa-1)r+1)}\le\sqrt{\kappa r/(\kappa-1)}$.
  It is nonincreasing in $k$, so it is at most $\sqrt{2r^2/(r+1)}<\sqrt{2r}$ for
  every $k\ge2r$. Choosing $\epsilon=1/k$
  costs only the change $r\mapsto r+1$ in the numerator.
\end{rema}

\begin{rema}[Relation to D-optimal design]\label{rem:dopt}
  The volume local search is precisely the Fedorov exchange local search for
  $0/1$ D-optimal design, as follows. In the notation of
  \cite[Section~1]{PFL2026}, D-Opt$(0/1)$ is
  \[
    \max\Big\{\ln\det\Big(\sum_{j\in N}x_jv_jv_j^\top\Big):\
    \sum_{j\in N}x_j=s,\ x\in\{0,1\}^n\Big\}
  \]
  for design points $v_1\,,\ldots,v_n\in\R^{d}$ and budget $s$ (we write $d$
  for the dimension, which \cite{PFL2026} calls $m$). Taking $d:=r$, $s:=k$, $v_j:=r_j$ (the columns of $R$, equivalently the rows of
  the design matrix $R^\top\in\R^{n\times r}$), and $x$ the indicator vector of
  $K$, the information matrix $\sum_jx_jv_jv_j^\top$ of the design
  (the Fisher information matrix of the corresponding linear model, up to
  the noise variance) is $R_KR_K^\top=G$, so the objective is
  $\ln\vol(K)$, finite exactly when $\rank R_K=r$. The Fedorov exchange moves
  $x\mapsto x+e_l-e_j$ with $x_j=1$ and $x_l=0$ are our swaps $K\mapsto
  K-j+l$; Lemma~\ref{lem:swap} is the rank-two update (Sherman--Morrison together
  with the matrix determinant lemma) that \cite[Section~4.4]{PFL2026} uses to
  evaluate them; and a $(1+\epsilon)$-local maximizer in the
  sense of Definition~\ref{def:loc} is a solution at which no exchange
  increases the D-Opt objective by more than $\ln(1+\epsilon)$. Our search is
  the best-improvement variant; the first-improvement variants of
  \cite[Section~4]{PFL2026} stop at points with the same property. The case
  $k=r$ is the \emph{saturated} design (as many design points as unknown
  parameters), and then the search is the determinant
  local search that \cite{XFLP2021,PFLX2025} use. We note two features of our setting. First, the choice of design points is canonical only up to
  an invertible linear map: replacing $R$ by $LR$ multiplies $\vol$ by
  $\det(L)^2$, so the design points may equally be taken as the columns of
  $A[S,\cdot\,]$ for any $S$ as in Proposition~\ref{prop:det} (rational data)
  or of $V_1^\top$ (where $\norm{V_1[j,\cdot\,]}_2^2$ is the leverage score
  used for initialization in Section~\ref{sec:exp}); the columns of $A$
  itself cannot be used directly when $m>r$, because then
  $\sum_{j\in K}a_ja_j^\top$ is singular. Second, only the $0/1$ version is
  relevant, because a repeated row of $H$ carries no information.

  Lemma~\ref{lem:lev} can be read as a local, discrete counterpart of the
  equivalence theorem of Kiefer and Wolfowitz \cite{KieferWolfowitz1960}: for a
  \emph{continuous} D-optimal design (nonnegative weights on the design
  points instead of a $0/1$ selection) of total weight $k$ with information
  matrix $B$, the \emph{variance function} $v\mapsto v^\top B^{-1}v$, which is
  proportional to the variance of the least-squares prediction at $v$,
  satisfies $\max_lr_l^\top B^{-1}r_l\le r/k$,
  while Lemma~\ref{lem:lev} gives $\alpha_l\le(r+\epsilon k)/(k-r+1)$ at every
  $(1+\epsilon)$-local optimum of the exchange search. For $\epsilon=0$,
  Lemma~\ref{lem:lev} is \cite[Lemma~9]{MSTX2019}, where it is the key step in
  their approximation analysis of local search for D-optimal design; our
  proof,
  through the leverages in $K+l$, is a minor variant of theirs, which sums
  the local-optimality inequalities over $i\in K$.

  Two consequences follow. First, any Fedorov-local optimum of D-Opt$(0/1)$
  with these design points and budget $k$, computed by any of the fast
  implementations, initializations, and exchange variants developed for
  D-optimal design (see \cite[Sections~4 and~5]{PFL2026} and the references
  there), yields an ah-symmetric reflexive generalized inverse with the
  guarantees of Theorem~\ref{thm:main} and Proposition~\ref{prop:spectral}.
  Second, exact methods for D-optimal design, such as the branch-and-bound
  algorithm of \cite{PFL2026}, can compute a global volume maximizer on
  moderate instances, which in particular is a $1$-local maximizer.
  Conversely, the proof of Proposition~\ref{prop:spectral} gives a statement
  about designs: at every $(1+\epsilon)$-local optimum $K$ of the exchange
  search with budget $k$, the information matrices of all points and of the
  selected points satisfy
  $\sum_{j\in N}r_jr_j^\top\preceq\big(1+(n-k)(r+\epsilon k)/(k-r+1)\big)
  \sum_{j\in K}r_jr_j^\top\,$. So Fedorov-local designs are, up to that
  factor, simultaneously good for every criterion that is monotone in the
  Loewner order (such as the A-criterion $\tr(B^{-1})$ and the E-criterion
  $\lambda_{\min}(B)$ of an information matrix $B$), relative to using all $n$
  points.
\end{rema}

The leverage bound of Lemma~\ref{lem:lev} also controls the spectral norm.
For $K$ with $\rank R_K=r$, let $H_K:=XU_1^\top\,$, where $X$ equals
$R_K^\dagger$ in the rows $K$ and is zero elsewhere. By
Lemma~\ref{lem:param}, $H_K$ is the ah-symmetric reflexive generalized
inverse of minimum Frobenius norm among those supported on $K$; for $|K|=r$
it coincides with the column-block solution.

\begin{prop}\label{prop:spectral}
  For every generalized inverse $H$ of $A$, $\norm{A^\dagger}_2\le\norm{H}_2\,$.
  For every $K$ with $\rank R_K=r$, writing $G:=R_KR_K^\top$ and
  $\Omega_K:=G^{-1/2}R_{N\setminus K}\,$, we have
  $\norm{H_K}_2\le(1+\norm{\Omega_K}_2^2)^{1/2}\,\norm{A^\dagger}_2\,$. If moreover
  $r\le k<n$, $\epsilon\ge0$, and $K$ is a $(1+\epsilon)$-local maximizer of
  $\vol$ over $k$-subsets, then $\norm{\Omega_K}_2^2\le(n-k)(r+\epsilon k)/(k-r+1)$.
\end{prop}

\begin{proof}
  For any generalized inverse $H$ we have $A^\dagger=A^\dagger AHAA^\dagger$,
  and $A^\dagger A$ and $AA^\dagger$ are orthogonal projectors; hence
  $\norm{A^\dagger}_2\le\norm{H}_2\,$. Because $U_1$ has orthonormal columns,
  $\norm{H_K}_2=\norm{R_K^\dagger}_2=\lambda_{\min}(G)^{-1/2}$, and
  $\norm{A^\dagger}_2=1/\sigma_r(A)=\lambda_{\min}(RR^\top)^{-1/2}$, where
  $\sigma_r(A)$ is the smallest positive singular value of $A$. Now
  \[
    RR^\top=G+R_{N\setminus K}R_{N\setminus K}^\top
    =G^{1/2}\big(I_r+\Omega_K\Omega_K^\top\big)G^{1/2}\preceq\big(1+\norm{\Omega_K}_2^2\big)\,G ,
  \]
  so $\lambda_{\min}(RR^\top)\le(1+\norm{\Omega_K}_2^2)\,\lambda_{\min}(G)$, which is
  the first bound. For the second,
  \[
    \norm{\Omega_K}_2^2\le\norm{\Omega_K}_F^2
    =\sum_{l\notin K}r_l^\top G^{-1}r_l=\sum_{l\notin K}\alpha_l\,,
  \]
  and Lemma~\ref{lem:lev} bounds each $\alpha_l\,$.
\end{proof}

\begin{rema}
  The first bound is an a~posteriori certificate, computable from $G$ and
  $R_{N\setminus K}$ alone. The second is the a~priori guarantee; at $k=r$
  it reads $\sqrt{1+r(n-r)(1+\epsilon)}\,\norm{A^\dagger}_2\,$, which is the
  bound that follows from \cite[Lemma~3.4]{Pan2000} for a local search on
  $|\det R_T|$ with tolerance $\sqrt{1+\epsilon}$. So
  Proposition~\ref{prop:spectral} recovers it without appeal to
  \cite{Pan2000}, and for $k>r$ the guarantee improves by roughly a factor
  $\sqrt{k-r+1}$. The step $\norm{\Omega_K}_2\le\norm{\Omega_K}_F$ can lose up to a
  factor $\sqrt{\min(r,n-k)}$, and in practice it does. On the instance S4 of
  Section~\ref{sec:exp}, with $\epsilon=0$ and $k/r=1$, $1.1$, $1.25$, and
  $1.5$, the a~priori bound on $\norm{H_K}_2/\norm{A^\dagger}_2$ is $100$,
  $28.6$, $17.0$, and $9.95$, the a~posteriori bound is $4.99$, $3.39$, $2.43$,
  and $1.65$, and the observed ratio is $3.55$, $2.26$, $1.59$, and $1.17$.
\end{rema}

\subsection{Polynomial running time}\label{sec:poly}

For rational data we measure input size by the number $\langle A\rangle$
of bits in a binary encoding of $A$. The argument follows the scheme of
\cite[Theorems~2.8 and~3.9]{XFLP2021} for $k=r$; the Cauchy--Binet formula
makes the bound on the number of swaps depend on $k$ only logarithmically\ for fixed $\epsilon$.

\begin{prop}\label{prop:poly}
  Let $A\in\mathbb{Q}^{m\times n}$ have rank $r$, let $r\le k<n$, and let
  $\epsilon\in(0,1]$ be rational. A $(1+\epsilon)$-local maximizer of $\vol$
  over $k$-subsets can be computed in time polynomial in $\langle A\rangle$
  and $1/\epsilon$. More precisely, there is an integer $\gamma\ge1$ with
  $\log_2\gamma\le 2\langle A\rangle$ such that, from a suitable starting
  $k$-subset, any sequence of improving swaps (each increasing $\vol$ by a
  factor larger than $1+\epsilon$) has length at most
  $(2r/\epsilon)\ln(ek\gamma^2)$.
\end{prop}

\begin{proof}
  \emph{Integral surrogate.} By Gaussian elimination, which by
  \cite[Section~3.3]{Schrijver1986} runs in polynomial time, compute a set $S$ of
  $r$ linearly independent rows of $A$ and a set $T_0$ of $r$ linearly
  independent columns. Let $B:=A[S,\cdot\,]\in\mathbb{Q}^{r\times n}$, let
  $d$ be the product of the denominators of the entries of $B$, and let
  $B':=dB$. Then $B'$ is an integer matrix whose entries are bounded in
  absolute value by some $\gamma$ with $\log_2\gamma\le2\langle A\rangle$. Define
  $\psi(K):=\det(B'_K{B'_K}^{\!\top})$. By Proposition~\ref{prop:det},
  $\psi(K)=d^{2r}\det(U_1[S,\cdot\,])^2\,\vol(K)$ for every $K$, with a
  positive factor independent of $K$. Hence $\psi$ and $\vol$ have the same
  $(1+\epsilon)$-local maximizers, and it suffices to run the local search
  on $\psi$. By the Cauchy--Binet formula,
  $\psi(K)=\sum_{T\subseteq K,|T|=r}\det(B'_T)^2$ is a nonnegative integer.

  \emph{Bounds on $\psi$.} Start from $K_0:=T_0\cup K'$, where $K'$ is any
  set of $k-r$ further columns. Because $\det B'_{T_0}$ is a nonzero integer,
  $\psi(K_0)\ge\det(B'_{T_0})^2\ge1$. For every $T$, Hadamard's inequality
  gives $\det(B'_T)^2\le\prod_{j\in T}\norm{B'_{\cdot j}}_2^2\le(r\gamma^2)^r$.
  Hence, using $\binom kr\le(ek/r)^r$, we get
  $\psi(K)\le\binom kr(r\gamma^2)^r\le(ek\gamma^2)^r$ for every $k$-subset $K$.

  \emph{Number of swaps.} Each accepted swap multiplies $\psi$ by more than
  $1+\epsilon$. Therefore, after $q$ swaps,
  $(1+\epsilon)^q<\psi_{\max}/\psi(K_0)\le(ek\gamma^2)^r$. Using
  $\ln(1+\epsilon)\ge\epsilon/2$ for $\epsilon\in(0,1]$, we get
  $q<(2r/\epsilon)\ln(ek\gamma^2)$.

  \emph{Cost of a swap.} At the current $K$, at most $k(n-k)$ candidate
  swaps are tested. Each test computes $\psi(K-j+l)$, the determinant of an
  $r\times r$ integer matrix with entries bounded by $k\gamma^2$, and compares it
  exactly with $(1+\epsilon)\psi(K)$. By \cite[Section~3.3]{Schrijver1986},
  this takes polynomial time. If no candidate passes the test, $K$
  is a $(1+\epsilon)$-local maximizer. Combining the three bounds proves the
  claim.
\end{proof}

\begin{coro}\label{cor:poly}
  For rational $A$ and $r\le k<n$, taking $\epsilon=1/k$ yields, in
  polynomial time, a $k$-subset $K$ with
  $z(K)\le\sqrt{k(r+1)/(k-r+1)}\;z^*$.
\end{coro}

\begin{rema}[Solving the restricted problem]\label{rem:socp}
  Problem~\eqref{eq:zK} can be written with rational data, without the
  singular-value decomposition. Let $T\subseteq K$ index $r$ linearly
  independent columns, let $H_T$ be the column-block solution over $T$, and
  let $Q\in\mathbb{Q}^{n\times(k-r)}$ have rows outside $K$ equal to zero and
  let the columns of $Q[K,\cdot\,]$ form a basis of the null space of
  $A[\,\cdot\,,K]$, which has dimension $k-r$.
  By Lemma~\ref{lem:param}, the ah-symmetric reflexive generalized inverses
  supported on $K$ are exactly the matrices $H=H_T+QYA[\,\cdot\,,T]^\top$ with
  $Y\in\R^{(k-r)\times r}$.
  Indeed, the matrices $X$ with $RX=I_r$ and zero rows outside $K$ are
  $X^{(T)}+QE$ with $E\in\R^{(k-r)\times r}$, where $X^{(T)}$ equals $R_T^{-1}$ in
  the rows $T$ and is zero elsewhere, and $X^{(T)}U_1^\top=H_T\,$. Moreover
  $EU_1^\top=YA[\,\cdot\,,T]^\top$ with $E=YR_T^\top\,$, because
  $A[\,\cdot\,,T]=U_1R_T\,$. The map $Y\mapsto QYA[\,\cdot\,,T]^\top$ is
  injective, so $z(K)=\min_Y\norm{H_T+QYA[\,\cdot\,,T]^\top}_{2{,}1}$ is an
  unconstrained second-order-cone program with rational data and a bounded
  set of optimal solutions. Writing each term $\norm{H_{i\cdot}}_2\le\tau_i$ with a
  second-order cone, as on \cite[p.~226]{NesterovNemirovskii1994}, turns it
  into a conic problem over a product of second-order cones. For such
  problems, \cite[\S6.2.3]{NesterovNemirovskii1994} establishes that
  interior-point methods reduce the inaccuracy of the current solution by
  an absolute constant factor in a number of iterations bounded by a
  polynomial in the total cone dimension, each iteration costing a
  polynomial number of arithmetic operations; \cite[\S6.1 and
  Theorem~3.2.1]{NesterovNemirovskii1994} gives the path-following analysis
  and the required a~priori bounds. In this sense,
  $z(K)$ can be approximated to any prescribed accuracy in polynomial
  time. In practice, the ADMM of \cite[Algorithm~2]{PFLX2025} applied
  to $R_K$ is considerably faster.
\end{rema}

\begin{rema}[Relation to volume sampling]\label{rem:vs}
  The analyses in \cite{AvronBoutsidis2013,DerezinskiWarmuth2017} of volume
  sampling (choosing a subset of columns with probability proportional to
  its squared volume) for column subset selection use the same leverage-sum
  step as Lemma~\ref{lem:lev}. A classical identity in the same circle of
  ideas, due to Berg \cite[Theorem~1 and Corollary]{Berg1986}, is that the
  Moore--Penrose solution of an arbitrary linear system, of any rank, is a
  convex combination of the solutions of its uniquely solvable $r\times r$
  subsystems, with weights proportional to the squared determinants of those
  subsystems. Ben-Tal and Teboulle \cite[Theorem~2.1]{BenTalTeboulle1990}
  treat the full-column-rank case, and \cite{DerezinskiWarmuth2017} gives the
  volume-sampling interpretation. In
  our setting this says that
  \[
    A^\dagger=\sum_{S,T}\frac{\det(A[S,T])^2}{\sum_{S',T'}\det(A[S',T'])^2}\;H_{S,T}\,,
  \]
  with both sums over the nonsingular $r\times r$ submatrices of $A$, and
  where $H_{S,T}$ is the block construction of \cite{FampaLee2018}, which
  places $A[S,T]^{-1}$ in rows $T$ and columns $S$. So $A^\dagger$ lies
  in the convex hull of the reflexive block generalized inverses.
  Grouping the terms by $T$, applying
  \cite[Theorem~2.1]{BenTalTeboulle1990} to each $A[\,\cdot\,,T]$, and
  using $\sum_S\det(A[S,T])^2=\det(A[\,\cdot\,,T]^\top A[\,\cdot\,,T])=\det(R_T)^2$
  (Cauchy--Binet) establishes, equivalently, that $A^\dagger$ is the
  $\det(R_T)^2$-weighted average of the column-block solutions $H_T\,$.
\end{rema}

\subsection{Limits of the certificate}\label{sec:limits}

Theorem~\ref{thm:main} was obtained by bounding the certificate $\rho(K)$
of Theorem~\ref{thm:cert}. The next result establishes that, for some instances,
\emph{every} $k$-subset has a large certificate. So no choice of $K$, and in
particular no refinement of the local search, can make the certificate
small in general.

The simplest such instances come from low coherence. For a matrix $A\in\mathbb{R}^{m\times n}$ with
unit columns $a_1\,,\ldots,a_n\,$, let
$\operatorname{coh}(A):=\max_{i\ne j}|\langle a_i\,,a_j\rangle|$ denote its
\emph{mutual coherence}; this is not the \emph{coherence} of
\cite{DMMW2012}, which is the largest leverage score. Also, considering $N$ partitioned into sets $B_1\,,\ldots,B_p$ such that, for each
$s$, $\{a_j:j\in B_s\}$ is an orthonormal basis of the range of $A$. Two such bases are \emph{mutually unbiased} if
 $|\langle a_i\,,a_j\rangle|=1/\sqrt r$ for all $i$ and $j$ in different bases. Unions of two or more mutually unbiased orthonormal bases of $\R^r$ have
coherence exactly $1/\sqrt r$. 

\begin{prop}\label{prop:coh}
  Let $A\in\R^{m\times n}$ have rank $r<n$ and unit columns $a_1\,,\ldots,a_n\,$,
  so that $\operatorname{coh}(A)>0$.
  Then every $K\subsetneq N$ with $\rank A[\,\cdot\,,K]=r$ has
  $\rho(K)\ge1/\operatorname{coh}(A)$.
\end{prop}

\begin{proof}
  Because $A=U_1R$ and $U_1$ has orthonormal columns,
  $\langle r_i\,,r_j\rangle=\langle a_i\,,a_j\rangle$ for all $i,j$; in
  particular each $r_j$ is a unit vector. Let $l\notin K$, and let $c$ satisfy
  $R_Kc=r_l\,$. Then
  \[
    1=\langle r_l\,,r_l\rangle=\sum_{i\in K}c_i\langle r_i\,,r_l\rangle
    \le\operatorname{coh}(A)\,\norm{c}_1\, ,
  \]
  so $\beta_K(l)\ge1/\operatorname{coh}(A)$.
\end{proof}

\begin{exam}\label{ex:mub}
  If $r$ is a power of $2$, let $\mathcal H_r$ be Sylvester's Hadamard
  matrix, defined by $\mathcal H_1:=[1]$ and
  $\mathcal H_{2t}:=\left[\begin{smallmatrix}\mathcal H_t&\mathcal H_t\\
  \mathcal H_t&-\mathcal H_t\end{smallmatrix}\right]$, and let
  $A:=[\,I_r\ \ r^{-1/2}\mathcal H_r\,]$; any Hadamard matrix of order $r$
  would do.
  Then $n=2r$, and $\operatorname{coh}(A)=1/\sqrt r$, so every $K$ with $r\le k<2r$ and
  $\rank A[\,\cdot\,,K]=r$ has $\rho(K)\ge\sqrt r$. If $r$ is a power of $4$,
  then, by \cite{CCKS1997}, there are $r/2+1$ mutually unbiased orthonormal
  bases of $\R^r$; taking all their vectors as columns gives
  $n=r(r/2+1)$ and $\operatorname{coh}(A)=1/\sqrt r$, so the same bound holds for every
  $k<r(r/2+1)$.
\end{exam}

Example~\ref{ex:mub} needs a Hadamard matrix of order $r$, which exists only
if $r\in\{1,2\}$ or $4\mid r$. Up to a factor $\sqrt2$, the discrete cosine
transform removes this restriction.

\begin{exam}\label{ex:dct}
  For $r\ge2$, let $C_r\in\R^{r\times r}$ be the orthonormal DCT-II matrix,
  with entries $(C_r)_{1j}:=1/\sqrt r$ and
  $(C_r)_{ij}:=\sqrt{2/r}\,\cos\big(\pi(2j-1)(i-1)/(2r)\big)$ for $i\ge2$,
  where $i,j\in\{1,\ldots,r\}$; by \cite{Strang1999}, it is orthogonal. Let $A:=[\,I_r\ \ C_r\,]$. Then $n=2r$, the columns of
  $A$ are unit vectors, columns within each block are orthogonal, and the
  inner products between the blocks are the entries of $C_r\,$. Hence
  $\operatorname{coh}(A)\le\sqrt{2/r}$, and Proposition~\ref{prop:coh}
  gives $\rho(K)\ge\sqrt{r/2}$ for every $K$ with $r\le k<2r$ and
  $\rank A[\,\cdot\,,K]=r$. Because every column of an orthogonal matrix is
  a unit vector, some entry has absolute value at least $1/\sqrt r$; so no
  instance $[\,I_r\ \ Q\,]$ with $Q$ orthogonal gives more than $\sqrt r$
  through Proposition~\ref{prop:coh}, and $C_r$ is within a factor $\sqrt2$
  of the best for every $r$.
\end{exam}

Example~\ref{ex:dct} does not extend to $k\ge2r$ in the same way: that needs
a third orthonormal basis whose mutual coherence with both $I_r$ and $C_r$
is $O(1/\sqrt r)$. The classical real transforms are not suitable: in our
computations, for each pair among the DCT-II, the DST-II, the DCT-IV, and
the Hartley basis, $\sqrt r$ times the mutual coherence grows with $r$,
from between $2.2$ and $2.4$ at $r=6$ to between $9.0$ and $10.0$ at
$r=100$, whereas each of them has mutual coherence at most $\sqrt{2/r}$
with $I_r\,$.

For $k\ge2r$ and general $r$, beyond the range of Examples~\ref{ex:mub} and~\ref{ex:dct}, we use
nets of the sphere instead. We use two classical facts, from Ball's notes \cite{Ball1997}, about the
unit sphere $\mathbb{S}^{r-1}\subset\R^r$. First, by \cite[Lemma~2.2]{Ball1997},
if $u$ is uniformly distributed on $\mathbb{S}^{r-1}$, then for every unit
vector $v$ and every $t\in[0,1)$, $\Pr(\langle v,u\rangle\ge t)\le e^{-rt^2/2}$.
Second, by \cite[Lecture~9]{Ball1997}, $\mathbb{S}^{r-1}$ has an $\eta$-net
with at most $(4/\eta)^r$ points; that is, every point of the sphere is within
Euclidean distance $\eta$ of some net point.

\begin{prop}\label{prop:lower}
  Let $r\ge2$ and $\eta\in(0,1)$, put $\theta:=1-\eta^2/2$, and let
  $A\in\R^{r\times n}$ have as its columns $a_1\,,\ldots,a_n$ the points of an
  $\eta$-net of $\mathbb{S}^{r-1}$, with $n\le(4/\eta)^r$. Then $\rank A=r$.
  Moreover, for every $k$ with $2k<\exp(\theta^2r/2)$ and every $K\subseteq N$
  with $|K|=k$ and $\rank A[\,\cdot\,,K]=r$, we have
  $\rho(K)\ge\theta\sqrt{r/(2\ln(2k))}$.
\end{prop}

\begin{proof}
  \emph{Rank, and reduction to the columns of $A$.} If all $a_j$ lay in a
  hyperplane with unit normal $u$, then $\norm{u-a_j}_2^2=2>\eta^2$ for all
  $j$, contradicting the net property at $u$. Hence $\rank A=r$. Because
  $m=r$, $U_1$ is square and orthogonal, so $R=U_1^\top A$. The quantity
  $\rho(K)$ is unchanged when $R$ is replaced by $LR$ for any nonsingular
  $L$, because $R_Kc=r_l$ if and only if $LR_Kc=Lr_l\,$. Hence we may compute
  $\rho(K)$ with the columns $a_j\,$.

  \emph{A direction nearly orthogonal to $K$.} Fix $K$, and fix
  $t\in\big(\sqrt{2\ln(2k)/r},\,\theta\big)$; this interval is nonempty by
  the assumption on $k$. For uniform $u\in\mathbb{S}^{r-1}$ and each $i\in K$,
  the event $|\langle a_i\,,u\rangle|\ge t$ is the union of two caps, so it
  has probability at most $2e^{-rt^2/2}$. By the union bound,
  $\Pr(\exists i\in K:|\langle a_i\,,u\rangle|\ge t)\le2ke^{-rt^2/2}<1$. So
  there is a unit vector $u$ with $|\langle a_i\,,u\rangle|<t$ for all
  $i\in K$.

  \emph{A column with large $\beta_K(l)$.} Pick a net point $a_l$ with
  $\norm{a_l-u}_2\le\eta$. Then
  $\langle a_l\,,u\rangle=1-\tfrac12\norm{a_l-u}_2^2\ge\theta>t$, so
  $l\notin K$. For any $c$ with $A[\,\cdot\,,K]c=a_l\,$,
  \[
    \theta\le\langle a_l\,,u\rangle=\sum_{i\in K}c_i\langle a_i\,,u\rangle
    \le t\,\norm{c}_1\,.
  \]
  Thus $\rho(K)\ge\beta_K(l)\ge\theta/t$. Letting $t$ decrease to
  $\sqrt{2\ln(2k)/r}$ gives the claim.
\end{proof}

Because $n$ does not enter the bound, $\eta$ can be taken as small as we
like, and then $\theta$ is as close to $1$ as we like.

The factor $\ln(2k)$ in Proposition~\ref{prop:lower} comes from the union
bound, which ignores that the $k$ events are bounded by the same random
direction. Replacing the union bound by \v{S}id\'ak's inequality, and the
uniform direction by a Gaussian one, replaces $\ln(2k)$ by a quantity of
order $\ln(k/r)$, with explicit constants. We use three facts about a
standard Gaussian vector $g\in\R^r$ and a standard normal variable $\zeta$:
\begin{enumerate}
\item by \v{S}id\'ak's inequality \cite{Sidak1967}, for symmetric slabs,
  $\Pr\big(\bigcap_{i\in K}\{|\langle a_i\,,g\rangle|\le
  s\}\big)\ge\prod_{i\in K}\Pr(|\langle a_i\,,g\rangle|\le s)$;
\item $\Pr(|\zeta|>s)\le\chi(s):=\min\{1,\sqrt{2/\pi}/s\}\,e^{-s^2/2}$ for
  $s>0$. Indeed, $\Delta(s):=e^{-s^2/2}-\Pr(|\zeta|>s)$ has $\Delta(0)=0$,
  $\Delta(s)\to0$ as $s\to\infty$, and $\Delta'(s)=e^{-s^2/2}\big(\sqrt{2/\pi}-s\big)$,
  so $\Delta$ first increases and then decreases, and hence $\Delta\ge0$; and
  $\Pr(\zeta>s)=\int_s^\infty(2\pi)^{-1/2}e^{-y^2/2}\,dy
  \le\int_s^\infty(y/s)(2\pi)^{-1/2}e^{-y^2/2}\,dy=(2\pi)^{-1/2}e^{-s^2/2}/s$,
  so $\Pr(|\zeta|>s)\le\sqrt{2/\pi}\,e^{-s^2/2}/s$ (the Mills-ratio bound);
\item by \cite[Lemma~1]{LaurentMassart2000},
  $\Pr(\norm{g}_2^2\le r-2\sqrt{r\xi})\le e^{-\xi}$ for $\xi>0$.
\end{enumerate}

\begin{prop}\label{prop:lowersharp}
  Let $A$, $\eta$, and $\theta$ be as in Proposition~\ref{prop:lower}, let
  $\delta\in(0,1)$ and $k\ge r$, and let $\bar s>0$ be the unique solution of
  $\chi(\bar s)=r\delta^2/(5k)$; note that $\bar s^2\le2\ln(5k/(r\delta^2))$.
  If $\bar\tau:=\bar s/\sqrt{(1-\delta)r}<\theta$, then every $K\subseteq N$
  with $|K|=k$ and $\rank A[\,\cdot\,,K]=r$ has $\rho(K)\ge\theta/\bar\tau$; in
  particular, $\rho(K)\ge\theta\sqrt{(1-\delta)r/(2\ln(5k/(r\delta^2)))}$.
\end{prop}

\begin{proof}
  The function $\chi$ is continuous and strictly decreasing on $(0,\infty)$,
  from $1$ to $0$, and $r\delta^2/(5k)\in(0,\frac15]$, so $\bar s$ exists and
  is unique; and $\chi(s)\le e^{-s^2/2}$ gives the bound on $\bar s^2$. Take
  $s>\bar s$ with $s/\sqrt{(1-\delta)r}<\theta$, and put $x:=\chi(s)$. Then
  $x<r\delta^2/(5k)\le\frac15$. Each $\langle a_i\,,g\rangle$ is standard
  normal, so by facts (1) and (2), and because $\ln(1-x)\ge-x/(1-x)$,
  \[
    \Pr\Big(\bigcap_{i\in K}\{|\langle a_i\,,g\rangle|\le s\}\Big)
    \ge(1-x)^k\ge\exp\Big(-\frac{kx}{1-x}\Big)
    >\exp\Big(-\frac{r\delta^2/5}{4/5}\Big)=\exp\big(-r\delta^2/4\big).
  \]
  By fact (3) with $\xi=r\delta^2/4$, the right-hand side is at least
  $\Pr(\norm{g}_2^2\le(1-\delta)r)$. Hence some $g$ has
  $|\langle a_i\,,g\rangle|\le s$ for all $i\in K$ and
  $\norm{g}_2^2>(1-\delta)r$. The unit vector $u:=g/\norm{g}_2$ then has
  $|\langle a_i\,,u\rangle|<t:=s/\sqrt{(1-\delta)r}<\theta$ for all
  $i\in K$, and the last part of the proof of Proposition~\ref{prop:lower}
  gives $\rho(K)\ge\theta/t$. Letting $s$ decrease to $\bar s$ gives the
  claim.
\end{proof}

For $k=\kappa r$, $\bar s$ depends only on $\kappa$ and $\delta$, so
Proposition~\ref{prop:lowersharp} gives $\rho(K)\ge\theta\sqrt{(1-\delta)r}/\bar s$,
which is of order $\sqrt r$ with no logarithm in $r$. This agrees, up to constants, with
what the volume estimate of Carl and Pajor
\cite[Remark~(5), p.~496]{CarlPajor1988} gives for the inradius of a
polytope with $2k$ vertices in the unit ball, but here all constants are
explicit. Proposition~\ref{prop:lower} is better for moderate $r$, and
Proposition~\ref{prop:lowersharp} for large $r$.

\begin{rema}\label{rem:tight}
  Take $k=\kappa r$ with $\kappa>1$. For a volume local maximizer with
  $\epsilon=0$, the proof of Theorem~\ref{thm:main} gives
  $\rho(K)\le\sqrt{\kappa r/(\kappa-1)}$. Write $\bar s_\kappa(\delta)$ for the
  solution of $\chi(\bar s)=\delta^2/(5\kappa)$, and let
  \[
    F_\kappa(r):=\sqrt{\frac{\kappa}{\kappa-1}}\cdot
    \min\Big\{\sqrt{2\ln(2\kappa r)}\,,\ \min_{0<\delta<1}
    \frac{\bar s_\kappa(\delta)}{\sqrt{1-\delta}}\Big\}.
  \]
  On the instances of Propositions~\ref{prop:lower}
  and~\ref{prop:lowersharp}, every $K$ has a certificate within the factor
  $F_\kappa(r)/\theta$ of this upper bound, whenever the corresponding
  proposition applies. So, as a bound on the certificate,
  Theorem~\ref{thm:main} cannot be improved in general by more than
  $F_\kappa(r)/\theta$, for any $\theta<1$; this is bounded independently of
  $r$. For $\kappa=2$, the inner minimum is about $3.16$, attained near
  $\delta\approx0.35$; $F_2(r)$ equals $3.84$ at $r=10$ and $4.19$ at
  $r=20$, never exceeds $4.47$, and the second term of the minimum is the
  smaller one once $r>36$. Taking $\theta=0.998$,
  Proposition~\ref{prop:lower} applies for every $r\ge7$ and
  Proposition~\ref{prop:lowersharp} for every $r\ge11$, so for $k=2r$ the
  factor is at most $4.5$ for every $r\ge7$. When $r$ is a power of $2$ and
  $1<\kappa<2$, or $r$ is a power of $4$ and $1<\kappa<r/2+1$,
  Example~\ref{ex:mub} gives every $K$ the certificate $\rho(K)\ge\sqrt r$,
  so the factor is at most $\sqrt{\kappa/(\kappa-1)}$, independent of $r$;
  for $\kappa=2$ it is $\sqrt2$. For every $r\ge2$ and $1<\kappa<2$,
  Example~\ref{ex:dct}, an instance with only $n=2r$ columns, gives every
  $K$ of rank $r$ with $|K|=k<2r$ the certificate $\rho(K)\ge\sqrt{r/2}$, so the factor
  is at most $\sqrt{2\kappa/(\kappa-1)}$, again independent of $r$; using
  the exact bound $\sqrt{kr/(k-r+1)}$ of Theorem~\ref{thm:main} instead of
  $\sqrt{\kappa r/(\kappa-1)}$, it is at most $\sqrt{2k/(k-r+1)}$. These results concern only the
  certificate, not the ratio $z(K)/z^*$ itself.
  Section~\ref{sec:bases} shows that the two can differ: on the
  instances of Example~\ref{ex:mub}, many supports with $\rho(K)=\sqrt r$
  have $z(K)=z^*$, while, for some of them and a suitable norm $\mu$ on the
  rows, the ratio $z_\mu(K)/z_\mu^*$ equals $\rho(K)=\sqrt r$.
\end{rema}

\subsection{Other norms on the rows}\label{sec:norms}

Nothing in Theorems~\ref{thm:cert} and~\ref{thm:main} is specific to the
Euclidean norm. Let $\mu$ be any norm on $\R^m$, and define
\begin{align*}
  z_\mu^*&:=\min\Big\{\sum_{i\in N}\mu(H_{i\cdot}):\ H\text{ an ah-symmetric
    reflexive generalized inverse of }A\Big\},\\
  z_\mu(K)&:=\min\Big\{\sum_{i\in K}\mu(H_{i\cdot}):\ H\text{ as above, with }
    H_{i\cdot}=0\text{ for }i\notin K\Big\}.
\end{align*}
For $\mu=\norm{\cdot}_2$ these are $z^*$ and $z(K)$. For
$\mu=\norm{\cdot}_1$ they are the minimum $1$-norms of ah-symmetric reflexive
generalized inverses, without and with the support restriction.

\begin{theo}\label{thm:norms}
  For every norm $\mu$ on $\R^m$ and every $K\subseteq N$ with
  $\rank R_K=r$, we have $z_\mu^*\le z_\mu(K)\le\rho(K)\,z_\mu^*\,$. In
  particular, the bound of Theorem~\ref{thm:main} and the polynomial-time
  guarantee of Corollary~\ref{cor:poly} hold for $z_\mu$ in place of $z$.
\end{theo}

\begin{proof}
  By Lemma~\ref{lem:param}, $H_{i\cdot}=X_{i\cdot}U_1^\top\,$, so
  $\mu(H_{i\cdot})=f(X_{i\cdot}^\top)$ with $f(x):=\mu(U_1x)$, a norm on
  $\R^r$ because $U_1$ has full column rank. Let $f^*$ be its dual norm. As in
  \eqref{eq:dual}--\eqref{eq:dualK},
  \[
    z_\mu^*=\max\{\tr\Lambda: f^*(\Lambda^\top r_j)\le1,\ j\in N\},\qquad
    z_\mu(K)=\max\{\tr\Lambda: f^*(\Lambda^\top r_j)\le1,\ j\in K\},
  \]
  by Lagrangian duality for the minimization of a finite convex function over
  a nonempty affine set. The proof of Theorem~\ref{thm:cert} uses only the
  triangle inequality, now for $f^*$, so it gives
  $z_\mu^*\ge z_\mu(K)/\rho(K)$. The remaining claims follow because
  Theorem~\ref{thm:main} and Corollary~\ref{cor:poly} were proved by bounding
  $\rho(K)$, which does not depend on $\mu$.
\end{proof}

\begin{rema}\label{rem:onenorm}
  For $\mu=\norm{\cdot}_1$ and $k=r$, Theorem~\ref{thm:norms} together with Corollary~\ref{cor:kr} gives the
  factor $r(1+\epsilon)$ of \cite[Theorem~3.7]{XFLP2021}, measured against
  the minimum $1$-norm over ah-symmetric reflexive generalized inverses.
  The result in \cite{XFLP2021} is stronger in one respect: it measures
  against the minimum over all ah-symmetric generalized inverses, without
  requiring reflexivity, that is, against the starting point of the
  rank--norm trade-off that \cite{FLP2021} studies. Whether that stronger
  benchmark can be used for $k>r$ is open. On the computational side, the restricted problem for
  $\mu=\norm{\cdot}_1$ is a linear program; with the rational
  parametrization of Remark~\ref{rem:socp}, \cite{Schrijver1986} establishes that
  it can be solved exactly in polynomial time. The lower bounds of
  Propositions~\ref{prop:coh}, \ref{prop:lower}, and~\ref{prop:lowersharp} also apply
  unchanged, because they concern only $\rho(K)$.
\end{rema}

\subsection{Matching the lower bound with D-optimal designs}\label{sec:dsupp}

The proof of Theorem~\ref{thm:main} represents $r_l$ through the
minimum-$\ell_2$ coefficients, which amounts to weighting the columns in $K$
uniformly. Other weights give other certificates.

\begin{prop}\label{prop:wcert}
  Let $K\subseteq N$ with $\rank R_K=r$, and let $w\in\R^K_+$ satisfy
  $\sum_{i\in K}w_i=1$ and $G_w:=\sum_{i\in K}w_ir_ir_i^\top$ nonsingular.
  Then $\rho(K)\le\max\{1,\max_{l\notin K}(r_l^\top G_w^{-1}r_l)^{1/2}\}$.
\end{prop}

\begin{proof}
  For $l\notin K$, put $c_i:=w_i\,r_i^\top G_w^{-1}r_l$ for $i\in K$. Then
  $R_Kc=G_wG_w^{-1}r_l=r_l\,$, and by the Cauchy--Schwarz inequality
  \[
    \norm{c}_1=\sum_{i\in K}\sqrt{w_i}\cdot\sqrt{w_i}\,|r_i^\top G_w^{-1}r_l|
    \le\Big(\sum_{i\in K}w_i\Big)^{1/2}\Big(\sum_{i\in K}w_i(r_i^\top G_w^{-1}r_l)^2\Big)^{1/2}
    =(r_l^\top G_w^{-1}r_l)^{1/2}. \qedhere
  \]
\end{proof}

With uniform weights, $G_w=G/k$, and the bound is
$\max_{l\notin K}(k\alpha_l)^{1/2}$, which is what the proof of
Theorem~\ref{thm:main} uses. The best weights come from D-optimal design. A \emph{design} on $N$ is a
vector $w\in\R^N_+$ with $\sum_{j\in N}w_j=1$, with information matrix
$G_w:=\sum_{j\in N}w_jr_jr_j^\top\,$, and a \emph{continuous D-optimal design}
is a maximizer $w^*$ of
$\ln\det(\sum_{j\in N}w_jr_jr_j^\top)$ over the simplex in $\R^N$. By the
Kiefer--Wolfowitz theorem \cite{KieferWolfowitz1960}, $w^*$ is D-optimal if
and only if $\max_{j\in N}r_j^\top G_{w^*}^{-1}r_j=r$. All D-optimal designs
have the same information matrix $G_{w^*}\,$, by strict concavity of
$\ln\det$, and this matrix lies in the convex hull of the matrices
$r_jr_j^\top$ in the space of symmetric $r\times r$ matrices, of dimension
$r(r+1)/2$. So, by Carath\'eodory's theorem, there is a D-optimal design with
at most $r(r+1)/2+1$ support points, and it can be obtained from any
D-optimal design by linear algebra.

\begin{theo}\label{thm:dsupp}
  Let $w$ be a design on $N$ with $G_w$ nonsingular and
  $\max_{j\in N}r_j^\top G_w^{-1}r_j\le Cr$ for some $C\ge1$; for $C=1$ these
  are exactly the continuous D-optimal designs. Then every
  $K\supseteq\supp(w)$ with $K\neq N$ satisfies $\rho(K)\le\sqrt{Cr}$, and
  hence $z_\mu(K)\le\sqrt{Cr}\;z_\mu^*$ for every norm $\mu$ on $\R^m$. A
  D-optimal design with $|\supp(w)|\le r(r+1)/2+1$ exists, and for every fixed
  $C>1$ a design with this support bound can be computed in a number of
  arithmetic operations polynomial in $n$ and $r$.
\end{theo}

\begin{proof}
  The bound on $\rho(K)$ is Proposition~\ref{prop:wcert} applied to the
  restriction of $w$ to $K$, and the bound on $z_\mu(K)$ then follows from
  Theorem~\ref{thm:norms}. Existence for $C=1$ was established above. For $C>1$,
  Khachiyan's algorithm \cite{Khachiyan1996} for the dual problem of finding
  a minimum-volume ellipsoid containing $\pm r_1\,,\ldots,\pm r_n$ computes
  such a $w$ in a polynomial number of arithmetic operations (in the real
  number model); its support can then be reduced to at most $r(r+1)/2+1$
  points without changing $G_w\,$, and hence without changing
  $\max_jr_j^\top G_w^{-1}r_j\,$, by Carath\'eodory's theorem.
\end{proof}

\begin{rema}\label{rem:dsupp}
  For $k=r(r+1)/2+1$ and $\epsilon=0$, Theorem~\ref{thm:main} already gives
  $\sqrt{kr/(k-r+1)}=\sqrt r\,(1+O(1/r))$, so Theorem~\ref{thm:dsupp} improves
  it only by a factor $1+O(1/r)$. Its point is exactness: on the instances of
  Example~\ref{ex:mub} with $r$ a power of $4$, where $n=r(r/2+1)>r(r+1)/2+1$,
  every $K\subsetneq N$ of rank $r$ has $\rho(K)\ge\sqrt r$, while a support of
  a D-optimal design gives $\rho(K)\le\sqrt r$. So, for supports of that
  size, $\sqrt r$ is exactly the best bound the certificate can give.
\end{rema}


\subsection{Certificate versus loss on unions of orthonormal bases}\label{sec:bases}

The lower bounds of Section~\ref{sec:limits} concern $\rho(K)$, and
Remark~\ref{rem:tight} stresses that they say nothing directly about
$z(K)/z^*$. The instances of Example~\ref{ex:mub} demonstrate that the gap between
the two can be as large as the certificate itself for the $2{,}1$-norm, and
yet that the certificate is attained exactly by the loss for another norm on
the rows. Throughout this subsection, the columns of $A$ are unit vectors,
and $N$ is partitioned into sets $B_1\,,\ldots,B_p$ such that, for each
$s$, $\{a_j:j\in B_s\}$ is an orthonormal basis of the range of $A$; in
particular $n=pr$. 

\begin{theo}\label{thm:bases}
  Let $A$ be as above, and let $K\subseteq N$ with $\rank R_K=r$.
  \begin{enumerate}
  \item\label{it:bases-a} $z^*=r$ and $z(K)\ge r$. If $K\supseteq B_s$ for some
    $s$, then $z(K)=z^*=r$, that is, the support $K$ has no $2{,}1$-norm
    loss.
  \item\label{it:bases-b} Suppose that $p\ge2$, that the bases are pairwise
    mutually unbiased, and that $K\supseteq B_s$ and $K\cap B_t=\emptyset$
    for some $s$ and $t$. Define the norm
    $\mu_K(y):=\max_{i\in K}|\langle a_i\,,y\rangle|+\norm{y-U_1U_1^\top y}_2$
    on $\R^m$. Then $\rho(K)=\sqrt r$, $z_{\mu_K}(K)=r$, and
    $z_{\mu_K}^*=\sqrt r$, so
    $z_{\mu_K}(K)/z_{\mu_K}^*=\rho(K)=\sqrt r$, while $z(K)/z^*=1$. The
    equality $\rho(K)=\sqrt r$ needs only that the bases are pairwise
    mutually unbiased, $p\ge2$, $K\supseteq B_s$ for some $s$, and $K\ne N$.
  \item\label{it:bases-c} If $k=qr$ with $1\le q\le p$, then every $K$ that is
    the union of $q$ of the sets $B_1\,,\ldots,B_p$ maximizes $\vol$ over all
    $k$-subsets of $N$; in particular, it is a $(1+\epsilon)$-local
    maximizer for every $\epsilon\ge0$.
  \end{enumerate}
\end{theo}

\begin{proof}
  As in the proof of Proposition~\ref{prop:coh},
  $\langle r_i\,,r_j\rangle=\langle a_i\,,a_j\rangle$ for all $i,j$, so each
  $r_j$ is a unit vector and each $Q_s:=R_{B_s}$ is an orthogonal $r\times r$
  matrix.

  (\ref{it:bases-a}) Let $X$ have $p^{-1}Q_s^\top$ in the rows $B_s\,$, for
  every $s$; then $RX=p^{-1}\sum_sQ_sQ_s^\top=I_r$ and each of its $n=pr$
  rows has norm $1/p$, so $z^*\le r$, while $\Lambda=I_r$ is feasible
  for~\eqref{eq:dual} and~\eqref{eq:dualK} because $\norm{r_j}_2=1$, so
  $z^*\ge r$ and $z(K)\ge r$; and if $K\supseteq B_s\,$, then $Q_s^\top$ in
  the rows $B_s$ and zeros elsewhere is feasible for~\eqref{eq:zK} with
  $2{,}1$-norm $r$.

  (\ref{it:bases-b}) The mutual coherence of $A$ is $1/\sqrt r$ and
  $K\ne N$, so Proposition~\ref{prop:coh} gives $\rho(K)\ge\sqrt r$.
  Conversely, for $l\notin K$, let $c$ have entries $c_i:=\langle
  r_i\,,r_l\rangle$ for $i\in B_s$ and $0$ for $i\in K\setminus B_s\,$. Then
  $R_Kc=Q_sQ_s^\top r_l=r_l\,$, and, because $l\notin B_s\,$, each
  $|c_i|=1/\sqrt r$, so $\norm{c}_1=\sqrt r$. Hence $\rho(K)=\sqrt r$.

  Next, $\mu_K$ is a seminorm, and if $\mu_K(y)=0$, then $y$ lies in the
  range of $U_1\,$, which is spanned by $\{a_i:i\in K\}$ because
  $\rank R_K=r$, and $y$ is orthogonal to every $a_i$ with $i\in K$; so
  $y=0$, and $\mu_K$ is a norm. As in the proof of Theorem~\ref{thm:norms},
  $f(x):=\mu_K(U_1x)=\max_{i\in K}|\langle r_i\,,x\rangle|$, because
  $a_i=U_1r_i$ and $U_1^\top U_1=I_r\,$. For $j\in K$ and every $x$ with
  $f(x)\le1$ we have $\langle r_j\,,x\rangle\le1$, so $f^*(r_j)\le1$; hence
  $\Lambda=I_r$ is feasible for the dual of $z_{\mu_K}(K)$, and
  $z_{\mu_K}(K)\ge r$. On the other hand, let $X$ have $Q_t^\top$ in the
  rows $B_t$ and zeros elsewhere. Then $RX=I_r\,$, and its row indexed by
  $j\in B_t$ is $r_j^\top\,$, with
  $f(r_j)=\max_{i\in K}|\langle r_i\,,r_j\rangle|=1/\sqrt r$, because
  $K\cap B_t=\emptyset$ and every $i\in K$ lies in a basis other than
  $B_t\,$. So $z_{\mu_K}^*\le r/\sqrt r=\sqrt r$. Together with
  Theorem~\ref{thm:norms},
  \[
    \sqrt r\le\frac{z_{\mu_K}(K)}{z_{\mu_K}^*}\le\rho(K)=\sqrt r,
  \]
  so equality holds throughout; moreover $z_{\mu_K}(K)\le r$, using
  $Q_s^\top$ in the rows $B_s$ and $f(r_j)=1$ for $j\in B_s\,$, so
  $z_{\mu_K}(K)=r$ and $z_{\mu_K}^*=\sqrt r$. Finally $z(K)/z^*=1$ by
  part~(\ref{it:bases-a}).

  (\ref{it:bases-c}) For every $k$-subset $K$ with $\rank R_K=r$, the
  matrix $G=R_KR_K^\top$ is positive definite with
  $\tr G=\sum_{j\in K}\norm{r_j}_2^2=k$, so, by the inequality between the
  arithmetic and geometric means of its eigenvalues,
  $\vol(K)=\det G\le(k/r)^r$. If $K$ is the union of $q$ of the sets $B_t\,$,
  then $G=\sum Q_tQ_t^\top=qI_r$ and $\vol(K)=q^r=(k/r)^r$.
\end{proof}

\begin{coro}\label{cor:bases}
  Let $r$ be a power of $4$, let $A$ be the union of the $r/2+1$ mutually
  unbiased orthonormal bases of $\R^r$ of Example~\ref{ex:mub}, and let
  $k=\kappa r$ with $\kappa\in\{1,\ldots,r/2\}$ (here $\kappa$ is an
  integer, because $K$ is a union of whole bases). Then the union $K$ of any
  $\kappa$ of the bases is a global maximizer of $\vol$ over $k$-subsets
  with $z(K)=z^*$, and, for the norm $\mu_K$ of
  Theorem~\ref{thm:bases}, $z_{\mu_K}(K)/z_{\mu_K}^*=\sqrt r$. For
  $\kappa\ge2$, the bound that Theorems~\ref{thm:main}
  and~\ref{thm:norms} give for this $K$ is
  $\sqrt{\kappa r^2/((\kappa-1)r+1)}=\sqrt{k/(k-r+1)}\,\sqrt r\le\sqrt{\kappa/(\kappa-1)}\,\sqrt r$,
  so, as a guarantee for all row norms, it is tight for the loss itself up
  to the factor $\sqrt{k/(k-r+1)}\le\sqrt{\kappa/(\kappa-1)}$, which is at
  most $\sqrt2$ for $\kappa=2$.
\end{coro}

\begin{proof}
  Because $K$ is the union of $\kappa\le r/2$ of the $r/2+1$ bases, it
  contains a basis and is disjoint from another one. Apply
  Theorem~\ref{thm:bases}.
\end{proof}

\begin{rema}\label{rem:basescomp}
  The norm $\mu_K$ depends on $K$; whether a single row norm can force a
  loss of order $\sqrt r$ on every volume local maximizer is not settled by
  Corollary~\ref{cor:bases}. For the $2{,}1$-norm, part~(\ref{it:bases-a})
  covers only supports that contain a basis, but computations suggest that
  the loss is usually small on these instances. For the Hadamard
  instances $[\,I_r\ \ r^{-1/2}\mathcal H_r\,]$ of Example~\ref{ex:mub} with
  $r\in\{4,8,16,32\}$ and $k\in\{\lfloor1.25r\rfloor,\lfloor1.5r\rfloor,
  \lfloor1.75r\rfloor,2r-1\}$, most volume local maximizers found from random
  starts had $z(K)/z^*=1$ to solver accuracy and $\rho(K)=\sqrt r$, but not
  all: for $r=16$ and $k=24$ we found a volume local maximizer with $12$
  columns from each basis, $\rho(K)=5$, and $z(K)/z^*=1.17$. At $k=r$, the global volume maximizers are the two
  bases, by Hadamard's inequality, with no loss, but the local search also
  stops at other supports, in fact from almost every random start  ($95$ of $100$ for $r=16$); in our
  runs their loss reached $2.11$ for $r=16$ and $1.87$ for $r=32$, and their
  certificates reached $11$ and $14.6$, respectively. For the three
  mutually unbiased bases of $\R^4$ and every $6$-subset of rank $4$, the
  largest $2{,}1$-norm loss is $1.26$, and $1.12$ over the volume local
  maximizers; for the local maximizers that contain one basis and miss
  another, the loss for $\mu_K$ equals $2=\sqrt r$, as
  Theorem~\ref{thm:bases} predicts. So, for the $2{,}1$-norm, lower bounds on
  the actual loss of a volume local maximizer with $k>r$ remain open, and
  the instances of Section~\ref{sec:limits} that are unions of bases do not
  seem able to provide them.
\end{rema}

\section{Numerical experiments}\label{sec:exp}

\medskip\noindent\textbf{Setup.} We use the test instances of
\cite{PFLX2025}, which are available together with its code at
\url{https://github.com/GabrielPonte/admmGinv}. They are dense random
matrices, generated by the procedure of \cite[Section~2.1]{FLPX2021}, with
$(m,n,r)=(100i,50i,25i)$ for the instances S$i$ and
$(m,n,r)=(1000i,500i,250i)$ for the instances L$i$, $i=1,\ldots,5$. All codes are in
Python and were run on the same machine, a single-core virtual
machine; times are wall-clock seconds and exclude the computation of the
singular-value decomposition, which all methods share.

The volume local search (``Vol'') maintains $\Gamma:=R^\top G^{-1}R$ by two
Sherman--Morrison updates per swap, so that every swap ratio of
Lemma~\ref{lem:swap} is read off $\Gamma$ and a best-improvement step costs
$O(n^2)$; $\Gamma$ is recomputed from scratch every $50$ swaps. Following
\cite[Section~4]{PFL2026}, it starts from $r$ linearly independent columns
found by pivoted QR, completed to $k$ columns by largest leverage score
$\norm{V_1[j,\cdot\,]}_2^2\,$, and stops when no swap improves $\vol$ by more
than a factor $1+10^{-10}$. The reported time of Vol includes the solution
of~\eqref{eq:zK} described next, and that of Vol$^+$ includes the time of
Vol. Problem~\eqref{eq:zK} on the resulting $K$ is then
solved by the ADMM of \cite[Algorithm~2]{PFLX2025} in the reduced form; on
L1--L5 the same ADMM applied to~\eqref{eq:zstar} reproduces the values of
$z^*$ and of $\norm{H^{2{,}1}_{opt}}_{2{,}0}$ reported in \cite{PFLX2025}, where
$H^{2{,}1}_{opt}$ denotes the computed minimizer of~\eqref{eq:zstar} and rows
of norm at most $10^{-5}$ are counted as zero. ``Vol$^+$''
polishes the output of Vol by dual-guided swaps: it tries replacing each of
the two rows of $X_K$ of smallest norm by each of the two columns with the
largest dual violation $\norm{\Lambda^\top r_l}_2\,$, re-solves~\eqref{eq:zK},
and accepts the best improving swap, for at most $10$ swaps or $200$\,s. The
quantity $\rho(K)$ is computed with one linear program per column (HiGHS).

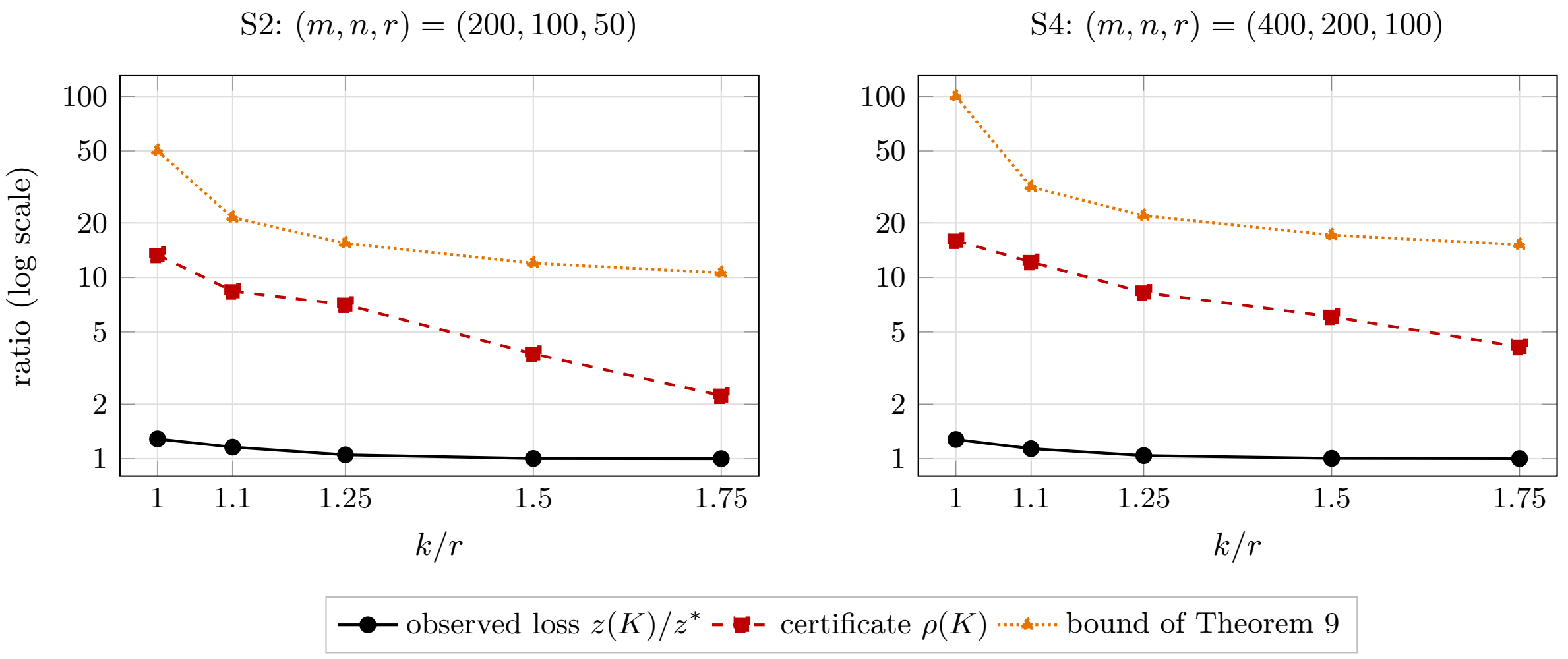
\begin{figure}[htbp]
  \centering
  \begin{tikzpicture}
    \begin{groupplot}[
      group style={group size=2 by 1, horizontal sep=1.6cm},
      width=0.47\textwidth, height=5.6cm,
      ymode=log, ymin=0.8, ymax=130,
      xmin=0.95, xmax=1.8, xtick={1,1.1,1.25,1.5,1.75},
      xticklabel style={/pgf/number format/fixed},
      xlabel={$k/r$}, grid=major, grid style={gray!25},
      tick label style={font=\footnotesize, black},
      label style={font=\small, black}, title style={font=\small, black},
      axis line style={black},
      legend style={font=\footnotesize, text=black, draw=gray!50,
        at={(-0.12,-0.3)}, anchor=north, legend columns=3},
      ytick={1,2,5,10,20,50,100}, yticklabels={1,2,5,10,20,50,100},
    ]
    \nextgroupplot[title={S2: $(m,n,r)=(200,100,50)$}, ylabel={ratio (log scale)}]
      \addplot[thick, black, mark=*] coordinates
        {(1,1.285) (1.1,1.158) (1.25,1.050) (1.5,1.002) (1.75,1.000)};
      \addplot[thick, red!75!black, dashed, mark=square*] coordinates
        {(1,13.35) (1.1,8.40) (1.25,7.09) (1.5,3.79) (1.75,2.23)};
      \addplot[thick, orange!90!black, densely dotted, mark=triangle*] coordinates
        {(1,50.00) (1.1,21.41) (1.25,15.44) (1.5,12.01) (1.75,10.62)};
    \nextgroupplot[title={S4: $(m,n,r)=(400,200,100)$}]
      \addplot[thick, black, mark=*] coordinates
        {(1,1.276) (1.1,1.136) (1.25,1.040) (1.5,1.004) (1.75,1.001)};
      \addlegendentry{observed loss $z(K)/z^*$}
      \addplot[thick, red!75!black, dashed, mark=square*] coordinates
        {(1,15.99) (1.1,12.19) (1.25,8.27) (1.5,6.08) (1.75,4.14)};
      \addlegendentry{certificate $\rho(K)$}
      \addplot[thick, orange!90!black, densely dotted, mark=triangle*] coordinates
        {(1,100.00) (1.1,31.62) (1.25,21.93) (1.5,17.15) (1.75,15.17)};
      \addlegendentry{bound of Theorem~\ref{thm:main}}
    \end{groupplot}
  \end{tikzpicture}
  \caption{Loss, certificate, and guarantee of the volume local
    maximizer against the number $k$ of allowed rows, for the instances S2
    and S4 (the data of Table~\ref{tab:k}; for S2 the values of $k$ are
    $50$, $55$, $62$, $75$, and $88$). The vertical axis is logarithmic.}
  \label{fig:k}
\end{figure}

\medskip\noindent\textbf{Loss versus $k$.} Table~\ref{tab:k} reports
$z(K)/z^*$ for the volume local maximizer on S2 and S4, together with the
certificate $\rho(K)$, the bound of Theorem~\ref{thm:main}, and the leverage
bound of Lemma~\ref{lem:lev}. Lemma~\ref{lem:lev} held in every run. The
loss decreases quickly: $25\%$ extra rows reduce it from about $1.28$ to
about $1.05$, and $50\%$ extra rows to below $1.005$. On S2 at $k=1.75r$,
$k$ exceeds the number of nonzero rows of $H^{2{,}1}_{opt}$ ($83$), and the
loss vanishes. Both $\rho(K)$ and the theoretical bound are loose, by factors
of about $2$--$13$ and $10$--$78$, respectively. Figure~\ref{fig:k} plots the same data
on a logarithmic scale. All three quantities decrease
with $k$, but at very different levels: the observed loss stays within $30\%$
of $1$ and is essentially $1$ from $k=1.5r$ on, the certificate is $2$ to
$13$ times larger, and the bound of Theorem~\ref{thm:main} lies a further
factor of $2$ to $6$ above the certificate.

\begin{table}[htbp]
  \centering
  \caption{Volume local search on $k$-subsets: instances S2,
    $(m,n,r)=(200,100,50)$, and S4, $(m,n,r)=(400,200,100)$.}\label{tab:k}
  \begin{tabular}{c|ccccc|ccccc}
    \hline
    & \multicolumn{5}{c|}{S2} & \multicolumn{5}{c}{S4}\\
    $k/r$ & $\frac{z(K)}{z^*}$ & $\rho(K)$ & Thm.~\ref{thm:main} & $\max\alpha_l$ & Lem.~\ref{lem:lev}
          & $\frac{z(K)}{z^*}$ & $\rho(K)$ & Thm.~\ref{thm:main} & $\max\alpha_l$ & Lem.~\ref{lem:lev}\\
    \hline
    1.00 & 1.285 & 13.35 & 50.00 & 6.32 & 50.00 & 1.276 & 15.99 & 100.00 & 5.12 & 100.00\\
    1.10 & 1.158 & 8.40  & 21.41 & 2.56 & 8.33  & 1.136 & 12.19 & 31.62  & 2.85 & 9.09\\
    1.25 & 1.050 & 7.09  & 15.44 & 1.48 & 3.85  & 1.040 & 8.27  & 21.93  & 1.27 & 3.85\\
    1.50 & 1.002 & 3.79  & 12.01 & 0.62 & 1.92  & 1.004 & 6.08  & 17.15  & 0.69 & 1.96\\
    1.75 & 1.000 & 2.23  & 10.62 & 0.20 & 1.28  & 1.001 & 4.14  & 15.17  & 0.36 & 1.32\\
    \hline
  \end{tabular}
\end{table}

\medskip\noindent\textbf{Comparison with ADMM heuristics.} For the target
number of rows $k=\lfloor\omega r+(1-\omega)\norm{H^{2{,}1}_{opt}}_{2{,}0}\rfloor$
used in \cite{PFLX2025}, Table~\ref{tab:h2h} compares the heuristics ADMM$_{2{,}0}$
and ADMM$_{2{,}1/0}$ of \cite[Algorithms~3 and~4]{PFLX2025} with Vol and
Vol$^+$. For
the two ADMM heuristics we use Python reimplementations of the code of
\cite{PFLX2025}, with its parameters, initialization, and placement of the
timer (in particular, as in \cite{PFLX2025}, the time of ADMM$_{2{,}1/0}$
excludes the run of ADMM$_{2{,}1}$ that initializes it, which took between
$0.5$ and $30$\,s); on every instance they reproduce the $2{,}1$-norms reported in
\cite[Table~2]{PFLX2025}, except where noted below. We imposed a time limit
of $1200$\,s on them; ADMM$_{2{,}1/0}$ reached it in five cases, in four of
which it had already met the row target, and its norm is within $0.12\%$ of
the value that \cite{PFLX2025} reports.

When few rows are removed ($\omega\le0.5$), Vol$^+$ and ADMM$_{2{,}1/0}$ are
comparable: Vol$^+$ is at least as good, and up to $0.68\%$ better, in seven of
the ten cases and worse by at most $0.12\%$ in the other three, and it is faster in
nine of the ten, by a factor of up to $4.7$. When the target is closer to
$r$ ($\omega\ge0.75$), Vol alone is better than ADMM$_{2{,}1/0}$ in every
case\ in which ADMM$_{2{,}1/0}$ met the row target, by $2.5\%$ to $27.3\%$\ (and by $34.6\%$ for L5 with $\omega=0.95$, where it did not), while being $68$ to $394$ times faster, and it
is better than ADMM$_{2{,}0}$ by $8.1\%$ to $37.9\%$; polishing rarely
changes it. In the most demanding case, L5 with $\omega=0.95$, ADMM$_{2{,}1/0}$ did
not reach the target within the time limit, whereas Vol needed $15$\,s.
Across all instances, Vol takes at most $17$\,s and Vol$^+$ at most $273$\,s.
The observed ratios $\norm{H}_{2{,}1}/z^*$ of Vol range from $1.006$ to
$1.50$, the certificate $\nu$ of Proposition~\ref{prop:post} from $1.52$ to
$9.67$, and the guarantee of Theorem~\ref{thm:main} from $26$ to $175$. The
guarantee explains why a few extra rows help, but it does not predict the observed
loss accurately.

\begin{table}[htbp]
  \centering
  \caption{On the instances L1--L5 of \cite{PFLX2025}: $2{,}1$-norms and
    times (seconds, same machine). ADMM$_{2{,}0}$ and ADMM$_{2{,}1/0}$ are the
    heuristics of \cite{PFLX2025} with a time limit of $1200$\,s; Vol is the
    volume local search followed by the restricted problem, and Vol$^+$
    adds dual-guided polishing. The last column is the certificate $\nu$ of
    Proposition~\ref{prop:post} for Vol. The smallest norm in each row is in
    bold.}\label{tab:h2h}
  \footnotesize\setlength{\tabcolsep}{3pt}
  \begin{tabular}{llr|rr|rr|rr|rr|r}
    \hline
    & & & \multicolumn{2}{c|}{ADMM$_{2{,}0}$} & \multicolumn{2}{c|}{ADMM$_{2{,}1/0}$}
      & \multicolumn{2}{c|}{Vol} & \multicolumn{2}{c|}{Vol$^+$} &\\
    inst.\ & $\omega$ & $k$ & norm & time & norm & time & norm & time & norm & time & $\nu$\\
    \hline
    L1 & 0.25 & 395 & 393.5 & 0.5 & 387.1 & 17.6 & 387.9 & 0.2 & \textbf{386.8} & 5.9 & 1.52\\
     & 0.50 & 347 & 411.6 & 0.8 & 398.7 & 10.0 & 399.6 & 0.2 & \textbf{396.0} & 4.6 & 2.23\\
     & 0.75 & 298 & 483.2 & 1.8 & 441.9 & 20.7 & 430.8 & 0.2 & \textbf{427.5} & 1.8 & 3.40\\
     & 0.80 & 288 & 512.1 & 3.0 & 474.1 & 20.5 & 442.0 & 0.3 & \textbf{440.4} & 1.2 & 3.37\\
     & 0.90 & 269 & 613.0 & 6.9 & 566.2 & 23.5 & \textbf{477.4} & 0.1 & \textbf{477.4} & 0.5 & 4.75\\
     & 0.95 & 259 & 729.5 & 17.0 & 643.8 & 21.9 & \textbf{517.1} & 0.1 & \textbf{517.1} & 0.5 & 7.15\\
    \hline
    L2 & 0.25 & 793 & 781.5 & 2.6 & \textbf{769.2} & 23.5 & 770.7 & 1.4 & \textbf{769.2} & 37.7 & 1.64\\
     & 0.50 & 695 & 807.8 & 5.2 & 790.0 & 115.1 & 790.2 & 1.0 & \textbf{787.5} & 25.2 & 2.11\\
     & 0.75 & 597 & 932.4 & 10.6 & 883.6 & 103.3 & \textbf{856.8} & 0.9 & \textbf{856.8} & 3.8 & 3.68\\
     & 0.80 & 578 & 1020.3 & 11.2 & 915.0 & 99.4 & \textbf{884.5} & 1.0 & \textbf{884.5} & 3.4 & 3.78\\
     & 0.90 & 539 & 1231.7 & 76.6 & 1138.1 & 152.0 & \textbf{985.4} & 0.8 & \textbf{985.4} & 3.0 & 5.66\\
     & 0.95 & 519 & 1734.5 & 266.0 & 1412.2 & 315.5 & \textbf{1077.4} & 0.8 & \textbf{1077.4} & 2.9 & 7.64\\
    \hline
    L3 & 0.25 & 1212 & 1184.4 & 9.0 & \textbf{1166.2} & 159.2 & 1170.1 & 3.6 & 1167.6 & 97.6 & 1.59\\
     & 0.50 & 1058 & 1249.9 & 18.2 & 1206.6 & 224.8 & 1205.5 & 3.3 & \textbf{1200.8} & 77.8 & 2.23\\
     & 0.75 & 904 & 1482.5 & 34.5 & 1395.7 & 353.9 & \textbf{1321.2} & 3.8 & \textbf{1321.2} & 11.9 & 3.27\\
     & 0.80 & 873 & 1609.5 & 45.6 & 1485.9 & 384.7 & \textbf{1366.2} & 4.2 & \textbf{1366.2} & 12.6 & 4.26\\
     & 0.90 & 811 & 2027.8 & 88.7 & 1826.9 & 510.8 & \textbf{1537.2} & 3.5 & \textbf{1537.2} & 10.7 & 5.77\\
     & 0.95 & 780 & 2580.9 & 516.4 & 2322.2 & 875.5 & \textbf{1687.6} & 3.8 & \textbf{1687.6} & 10.8 & 8.26\\
    \hline
    L4 & 0.25 & 1609 & 1582.1 & 19.0 & \textbf{1556.5} & 315.6 & 1559.5 & 8.1 & 1556.9 & 225.2 & 1.91\\
     & 0.50 & 1406 & 1665.2 & 47.3 & 1612.3 & 582.8 & 1612.7 & 6.4 & \textbf{1607.0} & 160.8 & 2.24\\
     & 0.75 & 1203 & 1967.3 & 76.9 & 1839.5 & 926.5 & 1786.8 & 7.6 & \textbf{1784.7} & 36.2 & 3.62\\
     & 0.80 & 1162 & 2117.3 & 87.5 & 1964.7 & 1008.5 & 1855.5 & 8.6 & \textbf{1852.7} & 35.3 & 4.09\\
     & 0.90 & 1081 & 2747.0 & 224.6 & 2417.1 & 1082.4 & \textbf{2080.3} & 7.3 & \textbf{2080.3} & 21.1 & 6.08\\
     & 0.95 & 1040 & 3592.4 & 436.2 & 3125.6 & 1200.1$^*$ & \textbf{2301.2} & 7.9 & \textbf{2301.2} & 21.6 & 9.67\\
    \hline
    L5 & 0.25 & 2035 & 1987.4 & 19.7 & \textbf{1958.6} & 335.6 & 1962.3 & 15.3 & 1960.3 & 272.9 & 1.82\\
     & 0.50 & 1773 & 2091.4 & 52.1 & 2033.1 & 1200.2$^*$ & 2038.0 & 17.0 & \textbf{2032.1} & 258.1 & 2.35\\
     & 0.75 & 1511 & 2478.9 & 106.2 & 2329.6 & 1159.2 & \textbf{2249.0} & 15.7 & \textbf{2249.0} & 48.4 & 3.69\\
     & 0.80 & 1459 & 2636.9 & 132.0 & 2461.2 & 1200.4$^*$ & \textbf{2335.0} & 14.7 & \textbf{2335.0} & 44.9 & 4.15\\
     & 0.90 & 1354 & 3563.9 & 297.8 & 3077.4 & 1200.2$^*$ & \textbf{2632.8} & 15.9 & \textbf{2632.8} & 43.1 & 6.48\\
     & 0.95 & 1302 & 4653.3 & 716.3 & 4475.0$^\dagger$ & 1200 & \textbf{2924.6} & 14.9 & \textbf{2924.6} & 40.2 & 8.75\\
    \hline
  \end{tabular}
  \par\smallskip
  \parbox{0.92\linewidth}{\footnotesize $^*$Time limit reached after the row
    target was met. $^\dagger$Row target not met within the time limit
    ($1373$ nonzero rows instead of $1302$).}
\end{table}

These experiments use a single family of random instances; other families,
and in particular instances whose $2{,}1$-norm optimum is row-sparse, may
behave differently.

\section{Concluding remarks}\label{sec:concl}

Allowing $k>r$ nonzero rows spans the range between the column-block
construction and $2{,}1$-norm minimization. A volume local search over
$k$-subsets, followed by a small convex problem, gives a guarantee that,
with $\epsilon=0$, is $r$ at $k=r$ and less than $\sqrt{2r}$ once $k\ge2r$. The guarantee is far above the losses we observed on the
instances of \cite{PFLX2025}: at most $1.005$ on S2 and S4 once $k\ge1.5r$,
and at most $1.047$ on L1--L5 at $\omega=0.5$, where $k\approx1.4r$. On those
instances the method, with a simple polishing step, is within $0.12\%$ of
the better ADMM heuristic of \cite{PFLX2025} when few rows are removed, and
gives $2{,}1$-norms $2.5\%$ to $27\%$ smaller, one to two orders of magnitude
faster, when the target is close to $r$. Our implementation borrows its
rank-one updates and its leverage-score initialization from the literature on
D-optimal design (Remark~\ref{rem:dopt}).

One gap in the theory deserves mention. For $k=\kappa r$ with
$\epsilon=0$, the certificate of a volume local maximizer is at most
$\sqrt{\kappa r/(\kappa-1)}$. When $r$ is a power of $4$ and
$\kappa<r/2+1$, Example~\ref{ex:mub} demonstrates that this is tight up to the
factor $\sqrt{\kappa/(\kappa-1)}$. For general $r$, however, our lower
bounds decrease like $\sqrt{r/\ln\kappa}$ as $\kappa$ grows
(Proposition~\ref{prop:lowersharp}), or apply only for $\kappa<2$ for every $r$ (Example~\ref{ex:dct}, and Example~\ref{ex:mub} with the
better constant when a Hadamard matrix of order $r$ exists). Whether instances with certificate of
order $\sqrt r$ for every $k$-subset exist for all $r$ and all $\kappa$,
ideally with polynomially many columns, is open; by a union bound over pairs and
Proposition~\ref{prop:coh}, $k+1$ random unit columns give order
$\sqrt{r/\ln k}$. For $\kappa<2$, Example~\ref{ex:dct} answers the question up to the factor $\sqrt2$ with $2r$ columns; for $\kappa\ge2$ and general $r$, one would need
a third orthonormal basis with mutual coherence $O(1/\sqrt r)$ to the first
two, and none of the classical real transforms we tried has this property.

There is also a gap between the certificate and the loss. 
Theorem~\ref{thm:bases} establishes that the lower bounds on $\rho(K)$ are
attained by the loss for a suitable norm on the rows, so that, when $r$ is a
power of $4$ and $k$ is a multiple of $r$ with $2r\le k\le r^2/2$, the
guarantee of Theorem~\ref{thm:norms} cannot be improved by more than the
factor $\sqrt{k/(k-r+1)}$ (Corollary~\ref{cor:bases}), but that they do not
transfer to the
$2{,}1$-norm: on unions of orthonormal bases every support containing a
basis has no $2{,}1$-norm loss. Whether the $2{,}1$-norm loss of a volume
local maximizer with $k\ge2r$ can grow with $r$, or is bounded by a
constant, is open.

\section*{Declaration on the use of AI}
During the preparation of this work, the authors used Claude to assist in
developing and checking all of the work. The authors take full responsibility for
the content of the article.

\bibliographystyle{plain+eid}
\bibliography{refs}

\end{document}